\newif\ifstandardtemplate\standardtemplatetrue%
\newif\ifelseviertemplate\elseviertemplatefalse%
\newif\ifspringertemplate\springertemplatefalse%
\newif\ifwileytemplate\wileytemplatefalse%
\newcommand{\mytitle}{A model and a finite element approximation of
  the mixed-dimensionality diffusion problem}
\newcommand{\myabstract}{%
  We present the formulation of a boundary value problem that models the coupled
  behavior of a three-dimensional diffusive solid with one-dimensional diffusive fibers embedded
  inside it. We introduce a variational statement of the problem that identifies the
  linked diffusive fields as energy minimizers under a coupling constraint. This saddle-point
  problem is proved to be well posed. Then, we introduce a finite element discretization
  of the proposed boundary value problem, and we prove the convergence of the finite element solution to
  the exact one. The most significant feature of this approximation is that the meshes
  of the bodies need not be conforming. Numerical examples confirm
  the theoretical results.}

\newcommand{\mypackages}{%
  \usepackage{amssymb}
  \usepackage{amsmath}
  \usepackage{amsthm}
  \usepackage{graphicx}
  \usepackage{natbib}
  \usepackage{todonotes}
  \usepackage{siunitx}
  \usepackage{booktabs}
  \usepackage{multirow}
  \usepackage{orcidlink}
  \usepackage[final]{changes} \definechangesauthor[name=Ignacio, color=orange]{IR}
  \graphicspath{{Figures/}{figures/}{./}}
  \usepackage{enumitem}
  \setlist[enumerate,1]{label = \emph{\alph*}),
                        ref   = \theenumi.\emph{\alph*}}
\usepackage{algorithm}
\usepackage{algpseudocode}}

\newcommand{\mymacros}{%
  \newcommand{\defined}{:=}
  \newcommand{\dev}{{\mathop{\mathrm{dev}}}}
  \renewcommand{\div}{{\mathop{\mathrm{div}}}}
  \newcommand{\mbs}[1]{\boldsymbol{##1}}
  \newcommand{\pairing}[2]{\langle{##1},{##2}\rangle}
  \newcommand{\dd}[2]{\frac{\mathrm{d} ##1}{\mathrm{d} ##2}}
  \newcommand{\pd}[2]{\frac{\partial{##1}}{\partial{##2}}}
  \newcommand{\set}[1]{\left\{##1\right\}}
  \newcommand{\trace}{{\mathop{\mathrm{tr}}}}
  \newcommand{\uptohere}{\centerline{\textcolor{blue}{\rule{6cm}{0.2cm}}}}
  \let\oldLambda=\Lambda\renewcommand{\Lambda}{\mathit{\oldLambda}}
  \let\oldGamma=\Gamma\renewcommand{\Gamma}{\mathit{\oldGamma}}
}
\ifstandardtemplate%
\documentclass[10pt,a4paper]{article}
\mypackages%
\usepackage[noblocks]{authblk}
\mymacros%

\theoremstyle{remark}
\newtheorem{theorem}{Theorem}
\newtheorem{remark}{Remark}

\title{\mytitle}
\author[1,2]{Ignacio Romero\orcidlink{0000-0003-0364-6969}}
\author[1]{David Portillo}

\affil[1]{Dept. of Mechanical Engineering, Universidad Politécnica de Madrid,
  Jos\'{e} Guti\'{e}rrez Abascal, 2, 28006 Madrid, Spain}
\affil[2]{IMDEA Materials Institute, Eric Kandel 2, 28096 Getafe, Madrid, Spain}

\begin{document}
\maketitle
\begin{abstract}
  \noindent\myabstract%
\end{abstract}
\fi


\section{Introduction}
\label{sec-intro}
It is not uncommon that bodies with widely different characteristic dimensions interact either
mechanically, thermally, chemically, or otherwise. For example, in mechanics, the stiffness and
strength of concrete civil structures is enhanced by using thin, long steel rebars. In physiology,
networks of very thin blood vessels transport blood into tissues and organs. Tree roots also cover
soil domains with network-like structures. Finally, advanced refrigeration systems in turbine blades
are based on narrow channels that deliver cooling air or fluids~\cite{heerden2022qv}. When analyzing these
problems, it proves convenient to employ continuum models for the medium and
one-dimensional reduced models for the embedded, slender bodies. This choice opens the door
to large computational savings when the problem is discretized, but naturally creates an unavoidable
problem: the coupling between the two types of bodies must be modeled, and then approximated.
Neither the formulation of well-posed problems of this type nor their approximation with
stable methods is simple.

Motivated by their relevance in multiple fields of application, several mixed-dimensionality models
have been proposed, and analyzed to show that they are well-posed. For example, by restricting the
solution spaces in the continuum and the embedded curve, the value of the concentration on the curve
and its average value on the cross section can be constrained to be equal \cite{angelo2008ca}.
Alternatively, the concentration on the curve can be extended to a cylindrical region, and
constrained to be equal to the concentration in the matrix restricted to the interface surface
\cite{laurino2019ps,kuchta2021on,berrone2022gi,berrone2023qp}. Other approaches might replace the
conducting curve with singular sources on the large body \cite{koch2020wx,gjerde2021uz,koch2022eu}. 

Similar problems arise in solid mechanics. To model the effects of thin fibers or inclusions
embedded in a deformable matrix, mixed dimensionality problems appear naturally. In this situation,
the presence of rotational degrees of freedom in the embedded structures complicates the
governing equations, although they are essentially equivalent to mixed dimensionality 
diffusion problems. The first \emph{models} for this kind of problem have been proposed by the authors
\cite{portillo2026ih}, including their stable discretization. Other numerical solutions for these mixed-dimensional solids have also been recently proposed
\cite{firmbach2023up,steinbrecher2022ge,steinbrecher2020jb,sky2024qd,hansbo2022fz}.

In this work we extend the ideas from our earlier work on embedded structures \cite{portillo2026ih}, adapting
them to diffusion problems. In this article, we use the language of thermal models, for concreteness, but the
results obtained can be used, almost without modification, for other diffusion problems governed by
Poisson's equation. We will show that stable formulations of mixed-dimensionality diffusive problems
can be obtained by constraining the unknown fields of the continuum and embedded conductor to be
equal in an interface volume. The key idea, originating from the work on the Arlequin method
\cite{dhia2001dh,dhia2005dh,quiao2011yx}, is to impose this constraint, not with standard Lagrange
multipliers, but with the natural (energy) inner product of the problem. We will show that this
critical modification is enough to ensure the well-posedness of the coupled problem. Moreover, this
desirable property carries over to Galerkin-type discretization almost directly.

An outline of the remaining sections of the article is the following. In Section~\ref{sec-model},
the mixed-dimensionality problem is described. A variational statement of the coupling is introduced
and the well-posedness of the problem is proven. The finite element discretization of this problem
is discussed in Section~\ref{sec-fem}. The convergence of the finite element solution to the exact
one is proven. Based on the proposed discretization, some illustrative examples of
mixed-dimensionality coupled problems are provided in Section~\ref{sec-examples}. The main results
of the article are summarized in Section~\ref{sec-conclusions}.

\section{A mixed-dimensionality thermal problem}
\label{sec-model}
In this section, we introduce the coupled problem that describes the thermal behavior of a 
slender conductor embedded into a three-dimensional continuum. Note
that modifying the latter with a two-dimensional domain would be
straightforward. Also, we note that the thermal problem is the prototypical example of diffusion,
and replacing it with a mass diffusion or an electrostatic problem is trivial.
In this article, for simplicity, we use the language of thermal transport. Also, other
elliptic problems arising in solid mechanics have somewhat different structure and we have
considered them elsewhere~\cite{portillo2026ih}.

\subsection{Geometry}
\label{subs-geo}
We start by presenting the geometry of the bodies that come into play in the problem of interest.
Without loss of generality, we will assume the simplest situation where a single long and thin
body is embedded in a three-dimensional body; the extension to multiple embeddings being trivial. To simplify the description, and avoid repeated periphrases,
we often refer to the three-dimensional body as the ``solid'' or the ``matrix'', and to the embedded slender body as the ``fiber''. 

The solid occupies a smooth bounded domain $\Omega\subset \mathbb{R}^3$ with boundary
$\partial\Omega$ and points denoted as $\mbs{x}$. The fiber is also a three-dimensional domain
$\mathcal{S}$ and we assume, for simplicity, that it can be generated by sweeping a circle $\Sigma$
along a smooth curve $\mathcal{C}$. If the length of $\mathcal{C}$ is $L$, points on this curve
can be assigned an arc-length coordinate $s\in[0,L]$ such that $\mbs{z}(s)$ is the one-parameter 
description of $\mathcal{C}$. To simplify the analysis, we will assume that the area of
$\Sigma$, denoted as $A$, is constant. The hypothesis that the fiber is slender is equivalent to
$A \ll |\Omega|^{2/3}$, with $|\Omega|$ denoting the total volume of the solid.

\subsection{The thermal problems}
\label{subs-thermal}
We are interested in studying the thermal behavior of the composite body consisting of a conductive
matrix and a conductive fiber embedded into it. When the two bodies are modeled as three-dimensional
entities, the coupled problem has a well-known mathematical structure: a Poisson equation with
possibly non-homogeneous conductivity can be used to find the thermal field everywhere. Moreover,
this boundary value problem is well-known and its discretization, using the finite element method or a similar one, follows without complication.

In this work, however, we are interested in modeling the fiber as a one-dimensional body, using the
equations of thermal conduction for curves and coupling its thermal field with the one in the
matrix. To make the notation clear and prepare the new results of Section~\ref{subs-coupling}, we
review next the governing equations of the two independent problems: the three-dimensional heat
conduction problem for the matrix and the one-dimensional counterpart for the fiber.

\paragraph{Three dimensional bodies.}
Starting from the thermal equilibrium of the matrix, let us consider a conductive continuum with
isotropic conductivity ${\kappa}>0$. The boundary $\partial\Omega$ of the matrix can be split into
two disjoint parts $\partial_D\Omega$ and $\partial_N\Omega$, where $\partial_D\Omega$ must
have nonzero measure. To proceed, let $L^2(\omega)$ be the space of functions that are
(Lebesgue) square-integrable over $\omega\subseteq\Omega$, for which the inner product
and norm are defined, respectively, as
\begin{equation}
  \label{eq-inner-l2}
  (u,v)_{L^2(\omega)} := \int_{\omega}
   u\,v \; dV\ ,
\quad
\|u\|_{L^2(\omega)} := (u,u)_{L^2(\omega)}^{1/2}\ .
\end{equation}
Similarly, let $H^1(\omega)$ be the Hilbert
space of functions in $L^2(\omega)$ with (weak) derivatives also in $L^2(\omega)$.
Functions $u,v$ in this space have inner product and associated norm defined,
respectively, as:
\begin{equation}
  \label{eq-inner-hilbert}
  (u,v)_{H^1(\omega)} := \int_{\omega}
  \left[
  u\,v + \ell^2\, \nabla u\cdot \nabla v
\right]\; dV\ ,
\quad
\|u\|_{H^1(\omega)} := (u,u)_{H^1(\omega)}^{1/2}\ ,
\end{equation}
where $\nabla$ denotes the gradient operator, the dot operation between the two vectors is just the
Euclidean scalar product, and $\ell$ is a characteristic length of
the integration domain~$\omega$. When $\omega\equiv\Omega$, for example, one can choose $\ell=|\Omega|^{1/3}$. 

If the body is subject to a volumetric heat supply $\bar{h}:\Omega\to \mathbb{R}$ and a
surface heat supply $\tilde{h}:\partial_N\Omega\to \mathbb{R}$, the equilibrium temperature on the
matrix is the field $u\in H^1_D(\Omega)$ that satisfies:
\begin{equation}
  \label{eq-body-minimizer}
  u = \arg\inf_{v\in H^1_D(\Omega)} I_{\Omega}[v]\ ,
\end{equation}
with $I_{\Omega}$ being the Dirichlet energy
\begin{equation}
  \label{eq-body-energy}
  I_{\Omega}[v] := \int_{\Omega} \frac{\kappa}{2} \| \nabla v(\mbs{x})\|^2 \; \mathrm{d} V
  - \int_{\Omega} \bar{h}(\mbs{x})\; v(\mbs{x}) \; \mathrm{d} V
  - \int_{\partial_N\Omega} \tilde{h}(\mbs{x})\; v(\mbs{x}) \; \mathrm{d} A\ .
\end{equation}
The set $H^1_D(\Omega)$ is the subset of $H^1(\Omega)$ of functions with vanishing trace on
$\partial_D\Omega$. Problem~\eqref{eq-body-minimizer} is well posed: the solution $u\in
H^1_D(\Omega)$ is unique and it depends continuously on the supplied heat (see, for example, \cite{hackbusch1992kl,evans1999wj}).

\paragraph{One dimensional bodies.}
In addition to modelling the temperature on the matrix, we would like to study the temperature on
the embedded fiber. For that, we start by recalling the formulation of the thermal problem in curves.
Once this problem is posed, we will discuss how the coupling between the continuum and the thin
inclusion is modeled.

As explained before, the fiber is a slender prismatic body $\mathcal{S}$ with constant cross
section~$A$. The centroids of the fiber are located on a one-parameter curve $\mathcal{C}$ of equation $\mbs{z}:[0,L]\to\mathbb{R}^3$ which, for simplicity, is
assumed to be smooth. The thermal conductivity per unit length of the fiber is $K>0$, assumed again for simplicity to be
constant. Let us note that the extension to nonsmooth curves or networks
thereof is straightforward as long as
the number of singular points has zero measure. Some heat per unit length
$\bar{Q}:\mathcal{C}\to\mathbb{R}$ is applied on the fiber and possibly some heat
$\tilde{Q}_\alpha$, with $\alpha=0$ and/or $\alpha=L$ might also be applied at the two ends of the curve.

To study the heat conduction on the fiber, let us define a second
Dirichlet functional $I_{\mathcal{C}}:H^1(\mathcal{C})\to \mathbb{R}$ of the form
\begin{equation}
  \label{eq-energy-curve}
  I_{\mathcal{C}}[\beta] := \int_{\mathcal{C}} \frac{K}{2} |\beta'(s)|^2 \; \mathrm{d}s
  - \int_{\mathcal{C}} \bar{Q}(s)\;\beta(s) \; \mathrm{d} s
  - \left[\tilde{Q}_{\alpha}\,\beta(\alpha)\right]_{\alpha=0}^L\ .
\end{equation}
To characterize the solution space of temperatures on the fiber we need to
define $L^2(\mathcal{C})$ and $H^1(\mathcal{C})$. These are, respectively,
the space of square-integrable functions on $\mathcal{C}$ and the space of functions
in $L^{2}(\mathcal{C})$ with derivatives also in $L^2(\mathcal{C})$. For future
reference, these two space have scalar products and norms that are parallel to
those defined in Eqs.~\eqref{eq-inner-l2} and ~\eqref{eq-inner-hilbert}. Namely,
for functions $\beta,\theta\in L^2(\mathcal{C})$, 
\begin{equation}
  \label{eq-inner-l22}
  (\beta,\theta)_{L^2(\mathcal{C})} := \int_{\mathcal{C}}
  \beta\,\theta \; ds\ ,
\quad
\|\beta\|_{L^2(\mathcal{C})} := (\beta,\beta)_{L^2(\mathcal{C})}^{1/2}\ ,
\end{equation}
and for functions $\beta,\theta\in H^1(\mathcal{C})$, 
\begin{equation}
  \label{eq-inner-hilbert2}
  (\beta,\theta)_{H^1(\mathcal{C})} := \int_{\mathcal{C}}
  \left[
  \beta\,\theta + R^2\,  \beta'\; \theta'
\right]\; ds\ ,
\quad
\|\beta\|_{H^1(\mathcal{C})} := (\beta,\beta)_{H^1(\mathcal{C})}^{1/2}\ ,
\end{equation}
with $R:=\sqrt{A/\pi}$.
To obtain the equilibrium temperature on the thermally conductive
curve we search for $\theta\in H^1(\mathcal{C})$ that solves
\begin{equation}
  \label{eq-curve}
  \theta = \arg\inf_{\beta\in H^1(\mathcal{C})} I_{\mathcal{C}}[\beta]\ .  
\end{equation}

\begin{remark}
Unless we add Dirichlet boundary conditions to $\theta$, problem~\eqref{eq-curve} has
a minimizer that is unique only up to a constant. One way to remove this non-uniqueness is to impose the value of
$\theta$ on, at least, one of the ends of the curve. However, we will leave the problem as it is
because we are interested in situations where the fiber is completely embedded inside the matrix $\Omega$
and the temperature at the ends of the curve is not given \emph{a priori}. We will later show that by
coupling appropriately the thermal fields $u$ and $\theta$, we will remove the non-uniqueness in
the temperature on the curve, bypassing the need for essential boundary conditions in its
formulation.
\end{remark}

\begin{remark}
In principle, there would be no problem to consider cases where
the fiber were not completely embedded
in the matrix. In these situations, it could be possible to impose the temperature at the end of
the fiber outside $\Omega$. The problem without any Dirichlet boundary conditions could be thought to be the hardest one to analyze, since its well posedness can only be shown if the
coupling with the body is strong enough to ensure uniqueness of solutions.  
\end{remark}

\subsection{Coupling the thermal fields}
\label{subs-coupling}
The variational formulations of thermal problem on subsets of $\mathbb{R}^3$ and three-dimensional
curves, summarized in Section~\ref{subs-thermal}, are standard. The first novelty of this article,
presented next, consists in modeling the \emph{coupling} between the thermal fields in
these two domains when the curve is embedded in the domain. This is non-trivial because the two
problems under consideration are formulated for different functional spaces, namely $H^1_D(\Omega)$
and $H^1(\mathcal{C})$. Moreover, we cannot simply impose that $u$ coincides with $\theta$ on
$\mathcal{C}$ because functions on $H^1_D(\Omega)$ do not have well defined traces on curves. To
follow this avenue, one would have to replace the Hilbert spaces employed before with weighted ones
where the projections of functions on the matrix onto functions on the curve are well defined
(see \cite{angelo2008ca}). 

Here, we follow a different strategy. Instead of attempting to \emph{project} the temperature $u$
onto the curve $\mathcal{C}$ with the hope of constraining it to be equal to $\theta$, we proceed in
the opposite direction. Following our recent work \cite{portillo2026ih}, we propose to \emph{lift}
the temperature $\theta$ from $\mathcal{C}$ to $\mathcal{S}$, the slender three-dimensional cylinder
that coincides with the fiber. Once this extension is built, it is reasonable to constrain that the
thermal fields in the intersection of two three-dimensional bodies be equal. In the past, other
groups have followed a similar approach \cite{laurino2019ps,kuchta2021on}, assuming a uniform
temperature distribution across each cross section of the fiber. The similarities with the present
work end there, since the variational formulations are different.

\begin{figure}[t]
  \centering
  \includegraphics[width=0.6\textwidth]{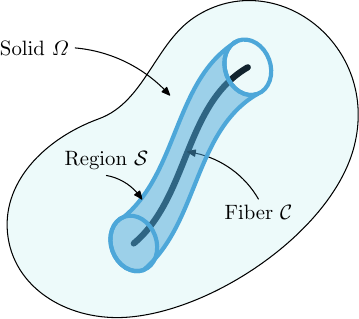}
  \caption{Solid matrix with an embedded thermal fiber. The fiber occupies a
    cylindrical volume that defines a curve of centroids.}
  \label{fig-buffer} 
\end{figure}

To construct the proposed lifting operation, consider the volume $\mathcal{S}$ that results from
sweeping a circle $\Sigma$ along the curve $\mathcal{C}$ (see Figure~\ref{fig-buffer}). At each point $\mbs{z}(s)\in \mathcal{C}$, the curve
$\mathcal{C}$ pierces the circle perpendicularly, and this circle has an area of value $A$.
Let $\Sigma(s)$ be the flat cross section of $\mathcal{S}$ that contains $\mbs{z}(s)$.
Points $\mbs{x}$ on this surface can be assigned
two coordinates $\xi^1,\xi^2$ by introducing a pair of unit directors $\mbs{d}_1(s),\mbs{d}_2(s)$
that span $\Sigma(s)$ and defining $\xi^{\alpha} = \mbs{d}_{\alpha}\cdot (\mbs{x}-\mbs{z}(s))$
for $\alpha=1,2$. Thus, there exist two projections $\Pi_{\mathcal{C}}:\mathcal{S}\to[0,L]$ and
$\Pi_{\Sigma}:\mathcal{S}\to \mathbb{R}^2$ such that, for all $\mbs{x}\in \Sigma(s)$
\begin{equation}
  \label{eq-projection}
  \Pi_{\mathcal{C}}(\mbs{x}) = s\ , \qquad
  \Pi_{\Sigma}(\mbs{x}) = (\xi^1,\xi^2)\ .
\end{equation}
We can combine these two projections and introduce a third projection $\Pi:\mathcal{S}\to
\mathcal{C}$ by 
\begin{equation}
  \label{eq-Pi}
  \Pi \mbs{x} = \mbs{z}(\Pi_{\mathcal{C}}(\mbs{x}))\ .
\end{equation}
Since all the points in $\Sigma(s)$ are projected by $\Pi$ onto $\mbs{z}(s)$, we can define the lifting operator
$\mathcal{L}:\mathcal{C}\to \mathcal{S}$ through the relationship
\begin{equation}
  \label{eq-lift}
  \Pi \circ \mathcal{L} = Id_{\mathcal{C}} ,
\end{equation}
where $Id_{\mathcal{C}}$ is the identity restricted to points on the fiber. The operator $\mathcal{L}$ is
set-valued and can be interpreted as the right inverse of $\Pi$. The left inverse of $\Pi$ is the
embedding operator $\mathcal{E}$ defined through the relationship
\begin{equation}
  \label{eq-left-inverse}
  \mathcal{E}\circ\Pi = Id_{\mathcal{S}}\ .
\end{equation}
It maps points on the product space $\mathcal{C}\times \Sigma$ onto the cylindrical region $\mathcal{S}$.

The projection and lifting operators can be used to study geometrical relations between
$\Omega,\mathcal{S}$ and $\mathcal{C}$. Also, they can be employed to lift fields defined on the curve
$\mathcal{C}$ to the region $\mathcal{S}$, giving to each section $\Sigma(s)$ the value at its
center $\mbs{z}(s)$. In fact, given an arbitrary field $\Phi\in
H^1(\mathcal{C})$ we can define its lift $\mathcal{L}\Phi$ to $\mathcal{S}$ by the composition
\begin{equation}
  \label{eq-lift-field}
  \mathcal{L} \Phi = \Phi\circ\Pi\ .
\end{equation}
Note that $\Pi$ is $C^{\infty}$ except on the boundary of $\mathcal{S}$, a region of zero measure.
Integrals of lifted fields can be evaluated on the curve, as in
\begin{equation}
  \label{eq-integral-lift}
  \int_{\mathcal{S}} \mathcal{L}\Phi(\mbs{x})\; \mathrm{d}V =
  \int_{\mathcal{C}} A\; \Phi(s)\; ds\ ,
  \qquad
  \int_{\mathcal{S}} \nabla(\mathcal{L}\Phi(\mbs{x}))\; \mathrm{d}V =
  \int_{\mathcal{C}} A\; \Phi'(s)\; \mbs{t}(s)\; ds\ ,
\end{equation}
where $\mbs{t}=\mbs{z}'$ is the unit tangent of the curve --- the second integrand is a vector,
since $\mathcal{L}\Phi$ varies only along $\mathcal{C}$ --- and where, to simplify the expressions,
we have selected $\ell\equiv R$ as the characteristic length of~$\mathcal{S}$.

Both identities in Eq.~\eqref{eq-integral-lift} are exact, and not only for a straight fiber. In a
cylindrical neighborhood of a smooth curve, the volume element is
$\mathrm{d}V = (1-k\,\xi)\,{d}A\,{d}s$, with $k$ the curvature of
$\mathcal{C}$ and $\xi$ the coordinate along the normal, while
$\nabla(\mathcal{L}\Phi) = \Phi'\,\mbs{t}/(1-k\,\xi)$. The factor $(1-k\,\xi)$ cancels in
the second identity, and integrates to $A$ over the symmetric cross section in the first, so no
approximation is involved as long as the tube does not self-intersect, i.e.\ $\kappa R<1$. The only
place where the curvature does enter is the gradient norm of Eq.~\eqref{eq-lift-norms} below, whose
cross-sectional factor is $\int_{\Sigma}(1-k\,\xi)^{-1}\mathrm{d}A = A\,(1+k^2\,R^2/4+\dots)$;
the first-order term vanishes by symmetry, so that identity holds up to a relative error
$O(k^2\,R^2)$, negligible for any slender fiber.

Suppose, finally, that we want to study the temperature of a body $\Omega$ that contains a fully
embedded conductive curve $\mathcal{C}$. The thermal equilibria of both bodies must be linked since
there might be heat exchange between them. Such a condition can thus be expressed by stating that the
coupled equilibrium will correspond to the thermal fields that minimize the thermal energies of the
body and curve, respectively, under the condition that the two fields be compatible on the
region~$\mathcal{S}$. The multiplier that enforces this compatibility is drawn from the space
\begin{equation}
  \label{eq-lambda-space}
  \Lambda := H^1(\mathcal{C})\ ,
  \qquad
  \| \mu \|_{\Lambda} := \| \mathcal{L}\mu \|_{H^1(\mathcal{S})}
  = \sqrt{A}\; \| \mu \|_{H^1(\mathcal{C})}\ ,
\end{equation}
that is, the space of thermal fields on the fiber, but equipped with the norm that the constraint
itself induces: a multiplier is only ever paired with fields lifted onto $\mathcal{S}$, so this is
its natural measure. The last identity in~\eqref{eq-lambda-space} is verified in
Eq.~\eqref{eq-lift-norms2} below, and shows that $\Lambda$ and $H^1(\mathcal{C})$ are the same set
of functions, normed differently by the constant factor~$\sqrt{A}$. We claim that the thermal
fields $u\in H^1_D(\Omega)$ and $\theta\in H^1(\mathcal{C})$ are
\begin{equation}
  \label{eq-coupled-problem}
  (u,\theta,\lambda) = \arg \inf_{v\in H^1_D(\Omega),\; \beta\in H^1(\mathcal{C})} \sup_{\mu\in \Lambda}
  L(v,\beta,\mu)
\end{equation}
with
\begin{equation}
  \label{eq-lagrangian}
  L(v,\beta,\mu) :=
      I_{\Omega}[v] +
      I_{\mathcal{C}}[\beta]
      +
      (\mathcal{L}\mu, \mathcal{L}\beta - v)_{H^1(\mathcal{S})}\ .
\end{equation}
Note that this functional incorporates the compatibility of the two temperature fields as a
constraint, in the $H^1$ sense, on the region~$\mathcal{S}$, where, for convenience,
we have selected $\ell \equiv R$.
We show next that this problem is well posed.

\subsection{Analysis}
\label{subs-analysis}
To prove that the saddle point problem~\eqref{eq-coupled-problem} is well-posed, let
us first collect some simple relationships and definitions. Let
$\phi\in H^1(\mathcal{C})$ and consider its lift $\mathcal{L}\phi$ to $H^1(\mathcal{S})$. Based on
the definition of $\mathcal{L}\phi$ we can easily verify that
\begin{equation}
  \label{eq-lift-norms}
  \| \mathcal{L}\phi \|_{L^2(\mathcal{S})} = \sqrt{A}\, \| \phi \|_{L^2(\mathcal{C})},
  \quad
  \| \nabla(\mathcal{L}\phi) \|_{L^2(\mathcal{S})} = \sqrt{A}\, \| \phi' \|_{L^2(\mathcal{C})}\ ,
\end{equation}
which together imply
\begin{equation}
  \label{eq-lift-norms2}
  \| \mathcal{L}\phi \|_{H^1(\mathcal{S})} = \sqrt{A}\, \| \phi \|_{H^1(\mathcal{C})}.
\end{equation}
Let us also define the space $U=H^1_D(\Omega)\times H^{1}(\mathcal{C})$, introduced to
collect the thermal fields on the domain and the embedded curve. The natural norm
on this product space is defined, for all $(v,\beta)\in U$, as
\begin{equation}
  \label{eq-u-norm}
  \| (v,\beta) \|_{U}^2 := \|v\|^2_{H^1(\Omega)} + A\; \|\beta\|^2_{H^1(\mathcal{C})}\ .
\end{equation}
With this notation we can now state the main result of this section.

\begin{theorem}
  \label{th-continuum}
Problem \eqref{eq-coupled-problem} is well posed.
\end{theorem}
\begin{proof}
To look for a triplet $(u,\theta,\lambda)\in H^1_D(\Omega)\times H^1(\mathcal{C})\times
\Lambda$ that can satisfy the stationarity conditions of the
Lagrangian~\eqref{eq-lagrangian}, let us first consider the stationarity conditions
of this functional. By taking functional derivatives of $L$, these conditions 
can be written as
  \begin{subequations}
    \label{eq-stationarity}
    \begin{align}
      a(u,\theta; v,\beta) + b(v,\beta; \lambda) &= \ell(v,\beta)\ , \label{eq-ab1} \\
      b(u,\theta; \mu) &= 0\ , \label{eq-ab2}
    \end{align}
  \end{subequations}
where $(v,\beta,\mu) \in H^1_D(\Omega)\times
H^1(\mathcal{C})\times \Lambda$ are arbitrary test functions.

In these equations we have introduced two bilinear forms, namely $a(\cdot;\cdot)$ and
$b(\cdot;\cdot)$, and one linear form $\ell(\cdot)$ defined, respectively, as
  \begin{equation}
    \label{eq-forms}
    \begin{aligned}
      a(u,\theta; v,\beta) &:=
                             \int_{\Omega} \kappa\, \nabla u\cdot\nabla v\; \mathrm{d}V
                             +
                             \int_{\mathcal{C}} K\, \theta'\,\beta'\; ds\ ,
      \\
      b(v,\beta; \lambda) &:= (\mathcal{L}\lambda,\mathcal{L}\beta-v)_{H^1(\mathcal{S})}\ ,
                            \\
      \ell(v,\beta) &:=
                      \int_{\Omega} \bar{h}\,v\;\mathrm{d}V + \int_{\partial_N\Omega} \tilde{h}\,v\; \mathrm{d}A
                      + \int_{\mathcal{C}} \bar{Q}\,\beta\; ds
                      + \left[\tilde{Q}_{\alpha}\,\beta(\alpha)\right]_{\alpha=0}^{L}\ .
    \end{aligned}
  \end{equation}
Problem~\eqref{eq-stationarity} falls within the class of \emph{mixed formulations} whose
well-posedness has been thoroughly studied~\cite{brezzi1991tn}. To ensure this property, we
must first recall that the kernel of the bilinear form $b(\cdot;\cdot)$ is the set
\begin{equation}
  \label{eq-kernel}
  \mathrm{ker}(b) =
  \left\{
    (v,\beta)\in U \ \hbox{such that}\
    b(v,\beta;\mu) = 0, \ \hbox{for all}\ \mu\in \Lambda
  \right\}.
\end{equation}
Three conditions are required to prove the well-posedness of the problem: (i) Both
$a(\cdot;\cdot)$ and $b(\cdot;\cdot)$ must be continuous, (ii) the bilinear form
$a(\cdot;\cdot)$ should be coercive in $\mathrm{ker}(b)$, and (iii) the bilinear form
$b(\cdot;\cdot)$ must satisfy the \emph{inf-sup} condition: there must exist a positive constant
$\gamma$ such that
\begin{equation}
  \label{eq-inf-sup}
  \inf_{\mu\in \Lambda} \sup_{(v,\beta)\in U}
  \frac{b(v,\beta;\mu)}{\|(v,\beta)\|_{U} \; \|\mu\|_{\Lambda}}
  \ge \gamma\ .
\end{equation}
To show the first condition, note that $a(\cdot;\cdot)$ is continuous on $U$ with a constant
proportional to $\max(\kappa,K/A)$, and, by Cauchy--Schwarz together with
Eqs.~\eqref{eq-lift-norms2} and~\eqref{eq-lambda-space}, we have that
\begin{equation}
  \label{eq-b-continuity}
  |b(v,\beta;\mu)| \le \| \mathcal{L}\mu \|_{H^1(\mathcal{S})}
  \left( \| \mathcal{L}\beta \|_{H^1(\mathcal{S})} + \| v \|_{H^1(\mathcal{S})} \right)
  \le \sqrt{2}\; \| \mu \|_{\Lambda}\; \|(v,\beta)\|_{U}\ .
\end{equation}
To prove the coercivity bound, let us first note that, by Poincaré's inequality, there exists a
constant $C_P>0$ such that, for all $u\in H^1_D(\Omega)$, 
  \begin{equation}
    \label{eq-coecivity-u}
    \int_{\Omega} |\nabla u |^{2} \;\mathrm{d}V \ge C_P\,\| u \|_{H^1(\Omega)}^2\ .
  \end{equation}
  Next, we observe that  if $(v,\beta)\in \mathrm{ker}(b)$ then, using
  the properties of the inner product as well as Eqs. \eqref{eq-lift-norms}
  and~\eqref{eq-lift-norms2}, we obtain
\begin{equation}
  \begin{aligned}
    0 &= b(v,\beta;-\beta) \\
      &= (\mathcal{L}\beta, v)_{H^1(\mathcal{S})} - (\mathcal{L}\beta,\mathcal{L}\beta)_{H^1(\mathcal{S})} \\
      &\le \| \mathcal{L}\beta \|_{H^1(\mathcal{S})} \, \| v \|_{H^1(\mathcal{S})} -
        \|\mathcal{L}\beta\|^{2}_{H^1(\mathcal{S})} \\
      &= \sqrt{A}\; \|\beta\|_{H^1(\mathcal{C})}\; \|v \|_{H^1(\mathcal{S})} -
        A\; \| \beta \|^2_{H^1(\mathcal{C})}\ .
  \end{aligned}
\end{equation}
Simplifying this relation we get
\begin{equation}
\| v \|_{H^1(\mathcal{S})} \ge \sqrt{A}\; \| \beta \|_{H^1(\mathcal{C})}\ .
\end{equation}
Using this bound, relation~\eqref{eq-coecivity-u}, and the fact that
$\mathcal{S}\subsetneq\Omega$, it follows that for all $(v,\beta)\in
\mathrm{ker}(b)$, there exists a constant $C>0$ such that
\begin{equation}
  \begin{aligned}
    a(v,\beta; v,\beta) &=
                          \kappa\; \| \nabla v \|^2_{L^2(\Omega)} + K \| \beta'\|^2_{L^2(\mathcal{C})} \\
                        &\ge \kappa\, C_P\, \| v \|^2_{H^1(\Omega)} \\
                        &\ge \frac{\kappa}{2} C_P\, \| v \|^2_{H^1(\Omega)} +
                          \frac{\kappa}{2} C_P\; A \|\beta\|^2_{H^1(\mathcal{C})} \\
                          &\ge C\;\kappa \|(v,\beta)\|^2_{U}\ ,
  \end{aligned}
\end{equation}
which proves the coercivity of $a(\cdot;\cdot)$ on the kernel of $b(\cdot;\cdot)$. To prove the
\emph{inf-sup} condition~\eqref{eq-inf-sup} it suffices to restrict the supremum to the particular
test pair $(v,\beta)=(0,\mu)$ and note that, by Eqs.~\eqref{eq-u-norm} and~\eqref{eq-lambda-space},
$\|(0,\mu)\|_{U} = \sqrt{A}\,\|\mu\|_{H^1(\mathcal{C})} = \|\mu\|_{\Lambda}$, so that
\begin{equation}
  \begin{aligned}
    \inf_{\mu\in \Lambda} \sup_{(v,\beta)\in U}
    \frac{b(v,\beta;\mu)}{\|(v,\beta)\|_{U} \; \|\mu\|_{\Lambda}}
    &\ge
    \inf_{\mu\in \Lambda}
      \frac{b(0,\mu;\mu)}{\|(0,\mu)\|_{U} \; \|\mu\|_{\Lambda}}
    \\
    &= \inf_{\mu\in \Lambda}
      \frac{\| \mathcal{L}\mu\|^2_{H^1(\mathcal{S})} }{\| \mu \|^2_{\Lambda}}
    \\
    &= 1 > 0\ .
  \end{aligned}
\end{equation}
With these results, the well-posedness of the coupled problem is now established.
\end{proof}

\begin{remark}
  \label{rk-scaling}
  The choice of norm in the multiplier space is not incidental. Had $\Lambda$ been normed with
  $\|\cdot\|_{H^1(\mathcal{C})}$ instead, the same argument would have given
  $\gamma=\sqrt{A}$, a constant with dimensions of length that vanishes in the slender limit
  $A\to0$, suggesting a loss of stability for thin fibers. No such loss occurs: the continuity
  constant in~\eqref{eq-b-continuity} would rescale by exactly the same factor, leaving the ratio
  that governs the stability estimates unchanged. The norm~\eqref{eq-lambda-space} simply makes
  this scaling explicit.
\end{remark}

\begin{remark}
  \label{rk-neumann}
  Theorem~\ref{th-continuum} assumes $|\partial_D\Omega|>0$, since the coercivity bound starts from
  Poincar\'e's inequality~\eqref{eq-coecivity-u} on $H^1_D(\Omega)$. The complementary situation is
  of practical interest and occurs in two of the examples of Section~\ref{sec-examples}: no
  essential data at all on the matrix, whose temperature is then determined \emph{only} through the
  tie to the fiber, while the fiber temperature is prescribed at one of its ends. Well-posedness
  survives, and the argument is worth recording because it is the coupling itself that supplies what
  the missing boundary condition would have.

  Let us therefore next assume that $\partial_D\Omega=\emptyset$, so that $v\in H^1(\Omega)$, and let the fiber
  temperature be prescribed at $s=0$, so that $\beta\in H^1_D(\mathcal{C}):=\{\beta\in H^1(\mathcal{C}),\;
  \beta(0)=0\}$.
  Since
  $\beta$ vanishes at one end, Poincar\'e's  inequality on the curve gives
  $\|\beta\|_{H^1(\mathcal{C})}\le C_{P}\,\|\beta'\|_{L^2(\mathcal{C})}$, so the fiber term
  of $a(\cdot;\cdot)$ alone controls $\|\beta\|_{H^1(\mathcal{C})}$. For the \emph{inf-sup} is
  to be satisfied, it suffices to select the multiplier space to be also $H^1_D$.
\end{remark}

\begin{remark}  We have
formulated problem~\eqref{eq-stationarity} with the goal of coupling the thermal behavior of the matrix and
the fiber, and finding thermal fields that model their joint behavior. It is important to
note, however, that the solution to this well-posed problem will not be, in general, identical to the true
solution of a three-dimensional slender body embedded inside another three-dimensional conductive body.

The coupled formulation admits, from the outset, two simplifications: first, the temperature
in the cross sections of the body $\mathcal{S}$ is constant, and second, there is an overlap
of the matrix and fiber, since the former is not removed when we define the latter. These two simplifications
are justified by the simplicity of the resulting model and its ability to reproduce the \emph{effective}
thermal behavior of the ensemble, but errors relative to the complete three-dimensional model should be
expected. This situation should not come as a surprise: all structural models
(beams, shells, plates, etc.) are also used to obtain, in a simplified fashion, approximate solutions to
the elasticity problem. As in the case of the coupled problem proposed in this section, the
merit of the structural models is to be judged by their simplicity and accuracy.
\end{remark}

\section{Finite element discretization}
\label{sec-fem}
In Section~\ref{sec-model}, we introduced the boundary value problem that models the thermal
behavior of coupled continua and curves, and proved that it is well posed. The
discretization of saddle point problems with finite elements is delicate, especially because the
crucial \emph{inf-sup} condition that is required to ensure the well-posedness of the continuous problem
need not be inherited by a Galerkin approximation~\cite{boffi2013jt}. In this section, we study a
finite element discretization of the constrained problem~\eqref{eq-coupled-problem} and prove that it is
stable and convergent. Moreover, in contrast with other mixed finite element formulations for which
the spaces of the primal and dual variables have to be carefully selected \cite{brezzi1991tn}, the
formulation introduced here is fairly robust in this respect. The solution spaces for the thermal
field in the matrix and the thermal field on the fiber can be arbitrary. The only restriction for
stability, as we will show, is that the space of multipliers be the same as the space of the
temperature on the fiber.

We start by defining the discretization spaces of functions defined on the solid $\Omega$ and
the curve $\mathcal{C}$. For that, we consider a mesh on each of these two bodies. The mesh
on~$\Omega$ partitions this volume into a set of volume elements $\mathcal{E}_{\Omega}=\{e_{i}\}$ connecting a finite
collection of nodes $\mathcal{N}_{\Omega}$ that define finite element shape functions
$\{N_a\}_{a\in \mathcal{N}_{\Omega}}$. The partition is assumed to be regular and we denote as
$h_{\Omega}$ a characteristic element dimension. Likewise, we partition the curve~$\mathcal{C}$ into
line elements $\mathcal{E}_{\mathcal{C}}=\{ E_i \}$ connecting the nodes
$\mathcal{N}_{\mathcal{C}}$, defining now finite element one-dimensional functions $\{M_a\}_{a\in\mathcal{N}_{\mathcal{C}}}$. We use the notation $h_{\mathcal{C}}$ to indicate the
characteristic element length in $\mathcal{E}_{\mathcal{C}}$. We stress that the volume and the
curve meshes are independently defined and, thus, possibly incompatible.

Next, we define finite element spaces of functions on $\Omega$ and $\mathcal{C}$. In the solid,
we introduce a set $V_h\subset H^1_D(\Omega)$ of finite element functions
\begin{equation}
  \label{eq-vh-continuum}
  V_h :=
  \left\{
  v_h(\mbs{x}) = \sum_{a\in \mathcal{N}_{\Omega}} N_a(\mbs{x})\; v_a
  \ \hbox{such that}\ v_h = 0 \ \hbox{on}\ \partial_D\Omega
  \right\}\ ,
\end{equation}
and on the curve we define
\begin{equation}
  \label{eq-wh}
  W_h :=
  \left\{
  \beta_h(s) = \sum_{a\in \mathcal{N}_{\mathcal{C}}} M_a(s)\; \beta_a
  \right\}\ .
\end{equation}
Finally, we introduce $U_h:= V_h\times W_h\subset U$.

Since $W_h$ is a subset of $H^1(\mathcal{C})$, the operator $\mathcal{L}$ defined in Eq.~\eqref{eq-lift}
can be used to lift functions from this finite element space to $H^1(\mathcal{S})$.

Using these objects, we now claim that the finite element solution of the mixed-dimensionality
diffusion problem are the temperature fields $(u_h,\theta_h)\in U_h$
and the multiplier $\lambda_h\in W_h$ that solve
\begin{equation}
  \label{eq-discrete-lagrangian}
  (u_h,\theta_h,\lambda_h) =
  \arg \inf_{(v_{h},\beta_h)\in U_h} \sup_{\mu_h\in W_h} L(v_h,\beta_h,\mu_h)
  \ .
\end{equation}

There is a key feature of the discrete problem that simplifies its analysis. The finite element
interpolation spaces for the temperature fields $(u_{h},\theta_h)$ are closed subspaces of the
functional spaces where the exact solution $(u,\theta)$ lie, namely, $U$. This is the rule for
(Bubnov-)Galerkin methods of one-field elliptic problems, but the exception for mixed finite elements.
For example, in the case of discrete formulations of Stokes' problem, discrete velocities are not
solenoidal, like the exact velocities. This situation complicates enormously the analysis of such
problems but not in the methods proposed here, as shown next.

\begin{theorem}
  \label{th-discrete}
  Problem~\eqref{eq-discrete-lagrangian} is well-posed and its solution converges to the solution of the continuum
problem \eqref{eq-coupled-problem}, i.e.,
\begin{equation}
  \label{eq-convergence}
  \lim_{h\to0} \| u_h - u \|_{H^1(\Omega)} = 0\ ,
  \qquad
  \lim_{h\to0} \| \theta_h - \theta \|_{H^1(\mathcal{C})} = 0\ .  
\end{equation}
\end{theorem}

\begin{proof}
  The proof of well-posedness is straightforward. The finite element
  problem~\eqref{eq-discrete-lagrangian} is of saddle point type, and thus analyzed using the same
  theory employed for the continuum problem in Theorem~\ref{th-continuum}. The variational equations of the problem are
  exactly~\eqref{eq-forms}, only now posed on the solution space $U_h=V_h\times W_h$ with
  multipliers in $W_h$. The steps of the proof of Theorem~\ref{th-continuum} can be replicated, just
  replacing the infinite dimensional spaces with their discrete counterparts, and the well-posedness
  of~\eqref{eq-discrete-lagrangian} follows. In particular, a discrete \emph{inf-sup} bound is
  satisfied: there exists a positive constant $\gamma_h$, independent of $h_{\mathcal{C}}$ and
  $h_{\Omega}$ such that
  \begin{equation}
    \label{eq-inf-suph}
    \inf_{\mu_h\in W_h} \sup_{(v_h,\beta_h)\in U_h}
    \frac{b(v_h,\beta_h;\mu_h)}{\|(v_h,\beta_h)\|_{U} \; \|\mu_h\|_{\Lambda}}
    \ge \gamma_h\ .
  \end{equation}
  In fact, the test pair used in the continuum proof, $(v_h,\beta_h)=(0,\mu_h)$, is available in the
  discrete setting as well, since $W_h\subset\Lambda$ implies $(0,\mu_h)\in U_h$ for every
  $\mu_h\in W_h$. The bound is therefore inherited from Eq.~\eqref{eq-inf-sup}, with $\gamma_h=1$ for every pair of
  meshes. This is precisely the point announced above: the multiplier space must coincide with the
  space of the temperature on the fiber, and no other compatibility between the two discretizations
  is required.

  The convergence result follows from the theory of mixed
  finite elements~\cite{brezzi1991tn}. Once the coercivity of $a(\cdot;\cdot)$ and the
  \emph{inf-sup} property of $b(\cdot;\cdot)$ are proven, the following estimate can be obtained:
  \begin{equation}
    \begin{split}
    \label{eq-error}
      \|(u_h-u,\theta_h-\theta)& \|_U + \|\lambda_{h} - \lambda\|_{\Lambda}
                                       \\
    \le &
    C \left(
      \inf_{(w_h,\beta_h)\in U_h}
      \|(w_h-u,\beta_h-\theta)\|_U
      +
      \inf_{\mu_h\in W_h}
      \|\mu_{h} - \lambda\|_{\Lambda}
      \right).
      \end{split}
  \end{equation}
  The convergence~\eqref{eq-convergence} is a consequence of Eq.~\eqref{eq-error} and the
  approximation properties of finite element spaces.  
\end{proof}

\subsection{Some remarks on the implementation}
\label{subs-implementation}
The finite element formulation described in the current section couples the thermal fields on the
matrix and fiber by appending a constraint to the Lagrangian~$L$ that enforces a (weak)
compatibility of temperature and temperature gradient in the cylindrical body~$\mathcal{S}$. In
practical implementations of the method, integrals over~$\mathcal{S}$ have to be performed with
sufficient accuracy. Note that at no point it is required that the
finite element meshes in the matrix and fiber be compatible.

To perform numerical quadrature over $\mathcal{S}$, one starts by selecting $N_{C}$ quadrature
points and weights $\{(s_i,W_i)\}$ on the fiber with $s_i\in[0,L]$ such that for every function $f:[0,L]\to\mathbb{R}$
\begin{equation}
  \label{eq-quadrature}
  \int_0^L f(s)\, ds \approx \sum_{i=1}^{N_{C}} f(s_i)\,W_{i}\ .
\end{equation}
Then, to approximate integrals on the cross section of $\mathcal{S}$, we introduce a second set of
quadrature points $\{(\xi_{1j},\xi_{2j},w_{j})\}_{j=1}^{N_A}$ such that, for every function defined
on the cross section $\Sigma$ we have
\begin{equation}
  \label{eq-quadrature2}
  \frac{1}{A}\int_{\Sigma} f(\xi^1,\xi^2)\, d\xi^1\, d\xi^2
  \approx
  \sum_{j=1}^{N_A} f(\xi^1_{j},\xi^2_{j})\, w_j\ ,
  \qquad
  \sum_{j=1}^{N_A} w_j = 1\ .
\end{equation}
A remark on this second rule is in order, because the coupling is only as good as it is. The
sampling has to represent the cross-sectional \emph{average} of the matrix field, since that is the
quantity the lifted fiber field is constrained against. The implementation described in
Appendix~\ref{appendix} places the sample points on two perpendicular diameters of $\Sigma$ rather
than over its area, a rule that is exact whenever $f$ varies at most linearly across the section ---
which is all the model itself resolves --- but not for its quadratic variation. As a result, the
solution of the discrete problem depends mildly on the transverse rule employed, through the second
moment it assigns to the cross section. This is a modelling choice, not an error, but it should be
kept in mind when the discrete solution is compared against a closed-form one, as it is in
Section~\ref{subs-quadratic-example}.
Combining these two quadrature rules, the integral of a lifted field $\mathcal{L}\Phi$, with $\Phi:\mathcal{C}\to\mathbb{R}$, can be
easily calculated as
\begin{equation}
  \int_{\mathcal{S}} \mathcal{L}\Phi(\mbs{x})\; dV = A \int_{\mathcal{C}} \Phi(s)\,ds
  \approx A \sum_i^{N_C} \Phi(s_i)\,W_i\ ,
\end{equation}
while the integral of an integrable function defined over the whole domain $\Omega$ but restricted
to $\mathcal{S}$ can be calculated as
\begin{equation}
  \int_{\mathcal{S}} u(\mbs{x})\,dV
  \approx
  \sum_{i=1}^{N_C} \sum_{j=1}^{N_A} (u\circ \mathcal{E})(\xi_{1j},\xi_{2j},s_i)\,W_i\,w_{j}\ .
\end{equation}
See Appendix~\ref{appendix} for further details on the implementation of the coupled problem.

\section{Numerical examples}
\label{sec-examples}
Finally, we illustrate the possibilities of the proposed
method by studying several mixed dimensionality examples.

\subsection{A consistency test}
\label{subs-patch-example}
First, we perform a simple consistency test on the formulation and its numerical implementation. An
exactly representable temperature field --- here a linear function --- is imposed on a body and an
embedded straight fiber. Since linear functions belong to the solution spaces in the matrix and
fiber, they must be reproduced to machine precision by the finite element method. Moreover, since
the field in the fiber is just the section of the three-dimensional temperature field, the Lagrange
multiplier must be identically zero everywhere. We use this same example to examine what happens to
that consistency as the fiber mesh is refined, and, also, to evaluate the key role played by the
gradient term in the coupling of Eq.~\eqref{eq-inner-hilbert}.

\begin{figure}[htbp]
  \centering
  \includegraphics[width=0.48\textwidth]{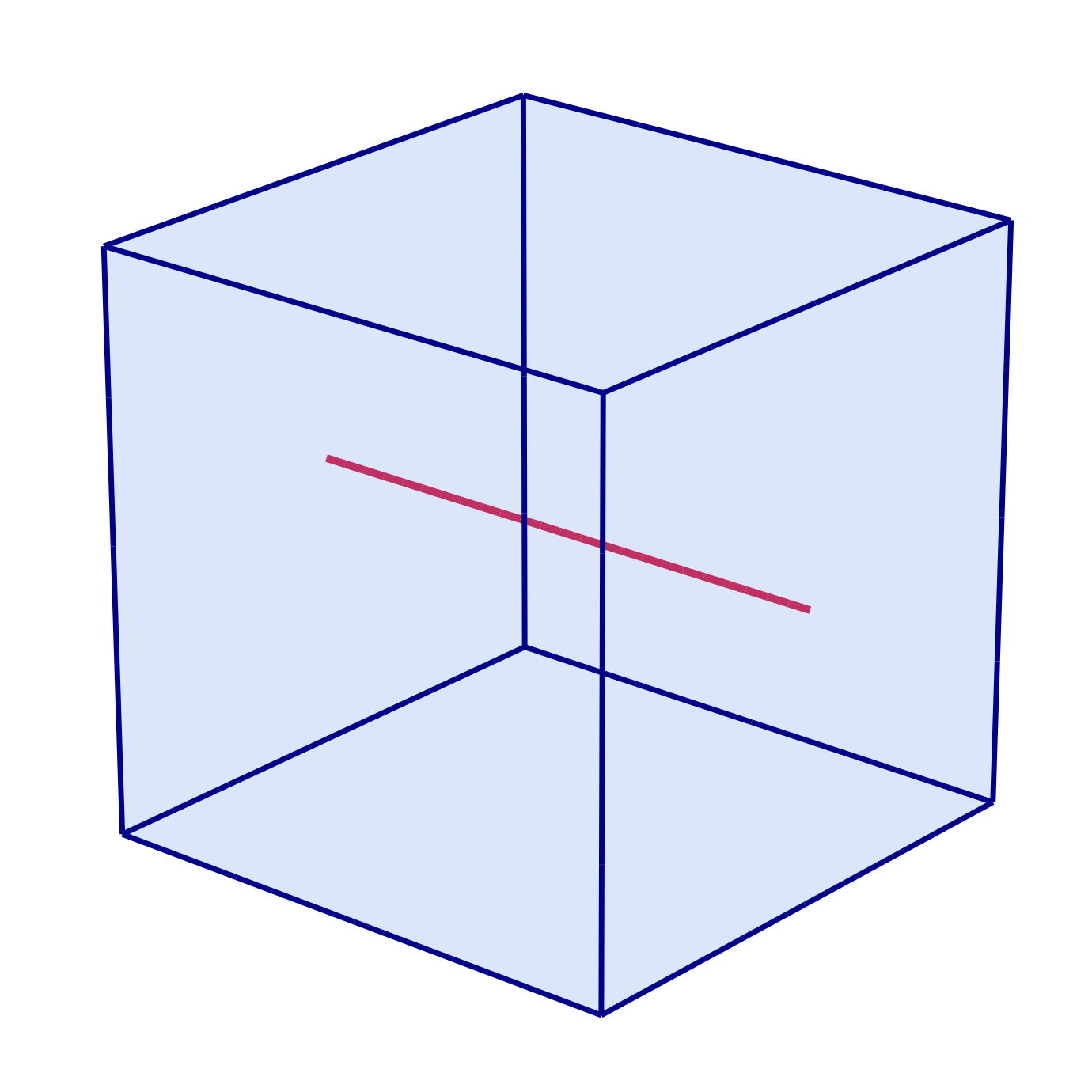}
  \includegraphics[width=0.48\textwidth]{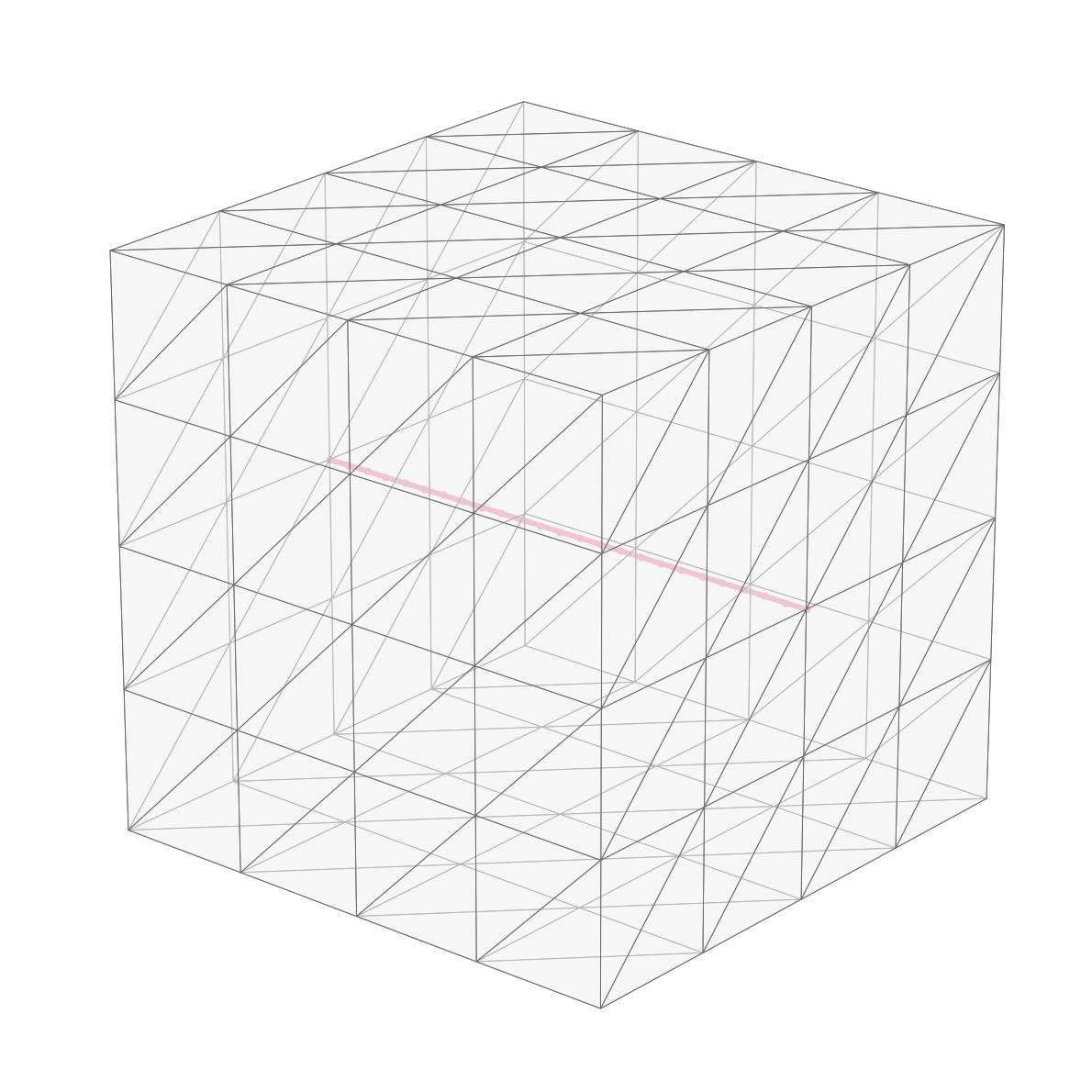}
  \caption{Geometry (left) and mesh (right, at a representative
    fiber refinement): straight fiber along the cube's axis, face center to face center. The
    solid mesh is held fixed throughout the refinement sweep described below; only the fiber
    mesh is refined.}
  \label{fig-patch-geometry} 
\end{figure}

As indicated, we study a thermally conducting cube of side $L=1$. A straight conductor is placed
from the center of one face of the cube to the opposite one. The embedded fiber has circular
cross section with radius $R=L/20$. The thermal conductivity of the solid is $\kappa=1$ and the
fiber linear conductivity is $K=20A\kappa$, with $A=\pi R^2$. A Cartesian coordinate system is
place at the center of the cube with axes parallel to the edges of the body. See Figure~\ref{fig-patch-geometry}.

A temperature $u_0=1$ is prescribed on the face $x=-L/2$ and on the fiber end that coincides
with it; likewise, $u_1=3$ is prescribed on the opposite face $x=+L/2$ and its fiber end;
the four lateral faces have zero thermal flux. The exact solution --- which coincides with the
finite element solution --- is the linear field
\begin{equation}
  u_h(x) = \theta_h(x) = u_0 + \frac{u_1-u_0}{L}\Bigl(x+\frac{L}{2}\Bigr) = 2+2x,
  \qquad \lambda_h =  0\ ,
  \label{eq-patch-exact}
\end{equation}
for every mesh. The temperature in the solid is a harmonic function with zero normal derivative on
the lateral faces, so it solves the solid problem. Since the thermal field in Eq.~\eqref{eq-patch-exact} has zero second derivative, it also
solves the fiber's own one-dimensional problem between the prescribed ends of the fiber, carrying a nonzero
axial flux $K\,(u_1-u_0)/L$. Because this field depends on $x$ alone and the cross sections
$\Sigma(s)$ are perpendicular to the $x$ axis, every coupling sample point on a given cross
section sees the same value of $u_h$, equal to $\theta_h$ there; the same holds for the axial
derivative. Thus, the $H^1(\mathcal{S})$ pairing
of $u_h-\theta_h$ with any multiplier vanishes identically. The solid and the fiber conduct \emph{in parallel},
while the tie between them must transmit \emph{exactly zero} heat; any nonzero $\lambda$ recovered
numerically is, therefore, discretization error.

One detail of the discretization deserves mention, because it is exactly the hypothesis on which
Theorem~\ref{th-discrete} rests. The temperature of the fiber is prescribed at its two ends, so the
two multipliers that would otherwise sit at those nodes have no matching temperature test function,
and the choice $(v_h,\beta_h)=(0,\mu_h)$ that proves the discrete \emph{inf-sup} bound is not
available for them. We therefore remove those two multipliers, which is to say we take $W_h$ to be
exactly the space of the discrete fiber temperature, as Section~\ref{sec-fem} requires. Doing so
changes nothing in this test --- the exact multiplier vanishes anyway --- but it is what makes the
constant reported below attain its theoretical value.

We solve this problem on a fixed solid mesh ($h_\Omega=0.25$) while refining the fiber
mesh, from $6$ to $384$ linear elements, comparing the default $H^1$ coupling of
Eq.~\eqref{eq-lagrangian} against a purely $L^2$ variant obtained by dropping the gradient
term of the constraint (equivalently, $\ell\equiv0$).

The quantity that governs the outcome is not the fiber element count, but rather the ratio
$R/h_{\mathcal{C}}$ between the coupling length $\ell\equiv R$ and the fiber element size. To
evaluate the stability of the formulation we compute the discrete \emph{inf-sup} constant
$\gamma_h$ of Eq.~\eqref{eq-inf-suph} algebraically \cite{chapelle1993wb,bathe2013tt}, as the square
root of the smallest eigenvalue of $P^{-1}\,C\,S^{-1}C^{T}$, where $C$ is the matrix of the
constraint and $S$ and $P$ are the Gram matrices of the norms $\|\cdot\|_U$ and $\|\cdot\|_{\Lambda}$
of Eqs.~\eqref{eq-u-norm} and~\eqref{eq-lambda-space}, respectively.
In addition, we calculate the maximum nodal value of the multiplier
$\lambda_h$.

\begin{table}[ht]
  \centering
  \small
  \begin{tabular}{lrrrrrrr}
    \toprule
    $R/h_{\mathcal{C}}$ & $0.30$ & $0.60$ & $1.20$ & $2.40$ & $4.80$ & $9.60$ & $19.20$ \\
    \midrule
    $\gamma_h$, $(\ell\equiv R)$ & $1.000$ & $1.000$ & $1.000$ & $1.000$ & $1.000$ & $1.000$ & $1.000$ \\
    $\gamma_h$, $(\ell\equiv0)$  & $5.3\cdot10^{-1}$ & $2.0\cdot10^{-1}$ & $5.5\cdot10^{-2}$ & $1.4\cdot10^{-2}$ & $3.6\cdot10^{-3}$ & $9.0\cdot10^{-4}$ & $2.3\cdot10^{-4}$ \\
    \midrule
    $\max|\lambda_h|$, $(\ell\equiv R)$ & $3.3\cdot10^{-13}$ & $6.5\cdot10^{-13}$ & $2.8\cdot10^{-13}$ & $3.8\cdot10^{-13}$ & $1.5\cdot10^{-12}$ & $1.1\cdot10^{-12}$ & $8.1\cdot10^{-13}$ \\
    $\max|\lambda_h|$, $(\ell\equiv0)$  & $4.8\cdot10^{-13}$ & $3.5\cdot10^{-12}$ & $3.7\cdot10^{-11}$ & $1.1\cdot10^{-10}$ & $3.7\cdot10^{-10}$ & $2.6\cdot10^{-9}$ & $9.1\cdot10^{-9}$ \\
    \bottomrule
  \end{tabular}
  \caption{Consistency test on the embedded fiber, solid mesh fixed
    ($h_\Omega=0.25$) and fiber mesh refined. $\gamma_h$ is the discrete \emph{inf-sup} constant of
    Eq.~\eqref{eq-inf-suph}; $\max|\lambda_h|$ is pure error, since the exact value is zero. The
    temperature is reproduced to $4\cdot10^{-11}$ at every level and in both formulations.}
  \label{tab-patch}
\end{table}

\begin{figure}[htbp]
  \centering
  \includegraphics[width=0.8\textwidth]{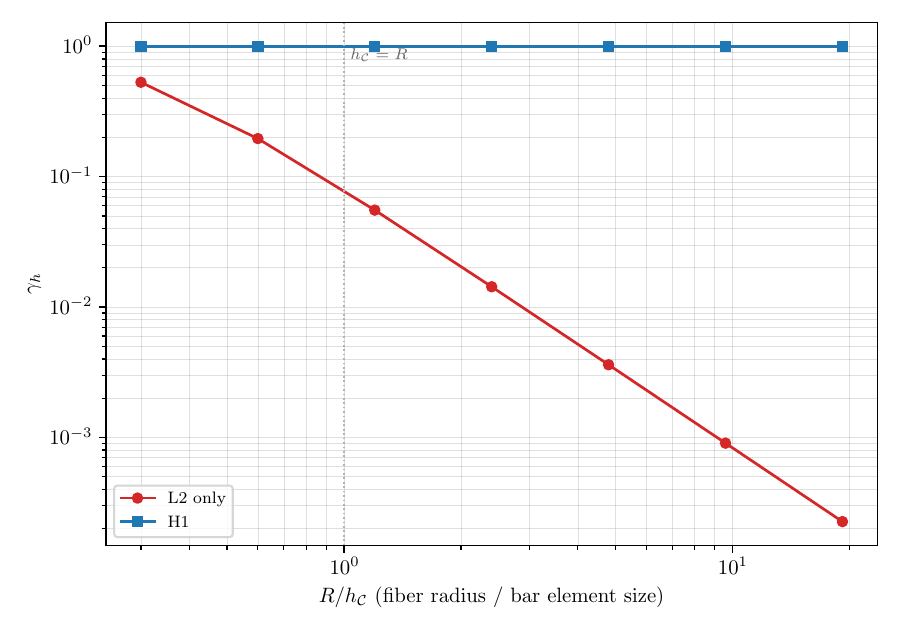}
  \caption{The discrete \emph{inf-sup} constant $\gamma_h$ of
    Eq.~\eqref{eq-inf-suph} across the fiber-refinement sweep. The proposed formulation gives
    $\gamma_h=1$ at every refinement, while the $L^2$ constrained formulation leads to
    $\gamma_h\to0$, pointing at an unstable numerical method.}
  \label{fig-patch-infsup}
\end{figure}

Table~\ref{tab-patch} and Figure~\ref{fig-patch-infsup} summarize the results of the analyses. The
proposed formulation passes the consistency test for all mesh refinements of the fiber: the
temperature is exact to round-off, the multiplier $\lambda_h$ never leaves the level of round-off,
and the \emph{inf-sup} constant satisfies $\gamma_h\approx 1$ at every one of the seven fiber
meshes, confirming Theorem~\ref{th-discrete}.

Instead, if the gradient term of the constraint is removed (i.e., $\ell\equiv0$), the discrete
\emph{inf-sup} constant decays quadratically under fiber refinement. Over the last four meshes the
computed values follow
\begin{equation}
  \label{eq-l2-decay}
  \gamma_h \simeq \frac{1}{12}\left(\frac{h_{\mathcal{C}}}{R}\right)^2
\end{equation}
to three significant digits, so the $L^2$-constrained method loses stability as soon as the fiber
mesh is refined below the fiber radius, and does so at a definite rate. The multiplier follows: with
the $H^1$ coupling $\max|\lambda_h|$ sits at $10^{-13}$--$10^{-12}$ irrespective of the mesh, while
with the $L^2$ coupling it grows by more than four orders of magnitude across the same sweep. Both
remain small in absolute terms only because the exact multiplier vanishes here, so what is being
amplified is round-off; the contrast between a flat sequence and a growing one is the meaningful
part, and Eq.~\eqref{eq-l2-decay} is what drives it.

\begin{figure}[htbp]
  \centering
  \includegraphics[width=0.48\textwidth]{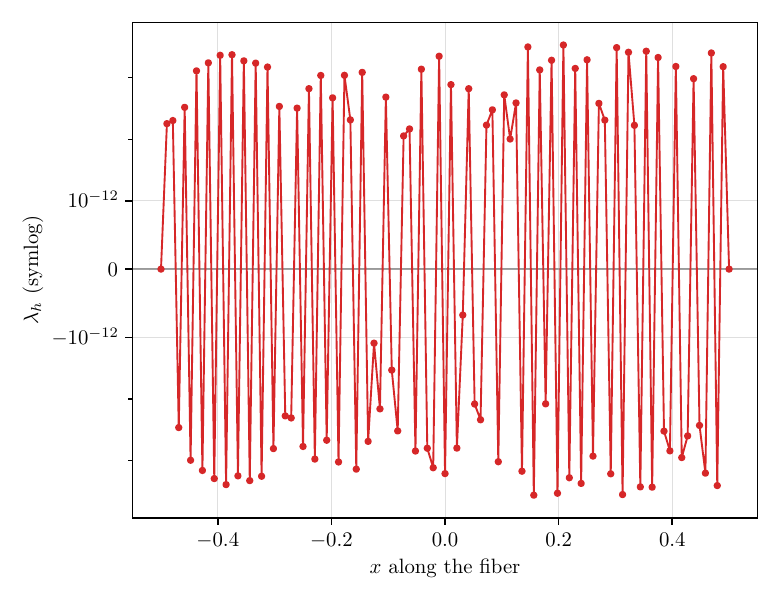}
  \includegraphics[width=0.48\textwidth]{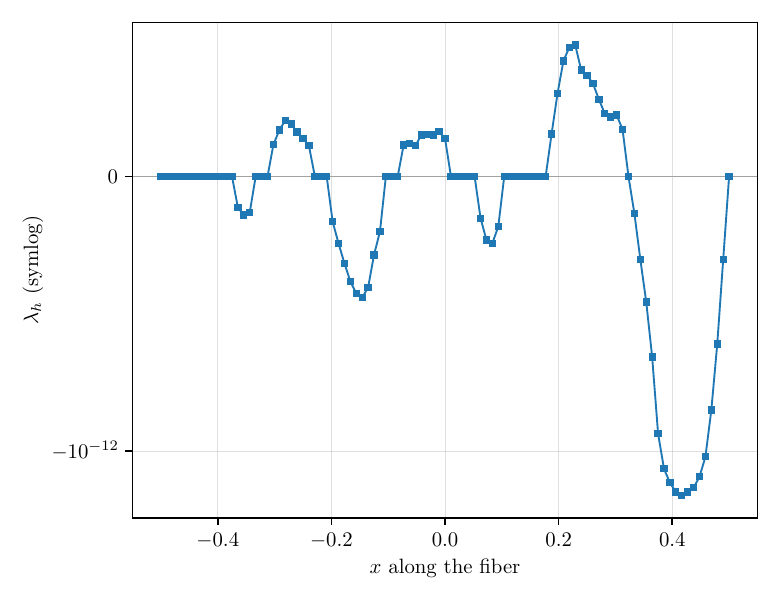}
  \caption{Value of the multiplier $\lambda_{h}$ along the
    fiber when using 96 elements in the fiber (symmetric logarithmic scale). $L^2$ formulation
    (left) and $H^1$ formulation (right).}
  \label{fig-patch-checkerboard}
\end{figure}

Figure~\ref{fig-patch-checkerboard} confirms that if no gradient is used in the
constraint ($\ell=0$), the multiplier oscillates along the fiber. In the
literature, this is often referred to as a \emph{checkerboard mode}. In the advocated formulation,
this \emph{mode} does not appear as a result of the unconditional stability.

\subsection{Convergence to an exact solution of the coupled problem}
\label{subs-quadratic-example}

The second example verifies the coupled discretization itself against a closed-form exact
solution --- of the coupled solid/fiber system, not of a classical single-body
idealization. Since the model introduced in Section~\ref{subs-coupling} enforces only a weak compatibility between the fiber and the matrix rather than a Dirichlet condition
on an actual cylindrical surface, there is, in general, no reason for the coupled solution to
coincide with the solution of any classical problem, even in the limit of mesh
refinement; the only rigorous way to verify the discretization is therefore against an exact solution
of the coupled formulation itself, coupling operator included.

\begin{figure}[htbp]
  \centering
  \includegraphics[width=0.31\textwidth]{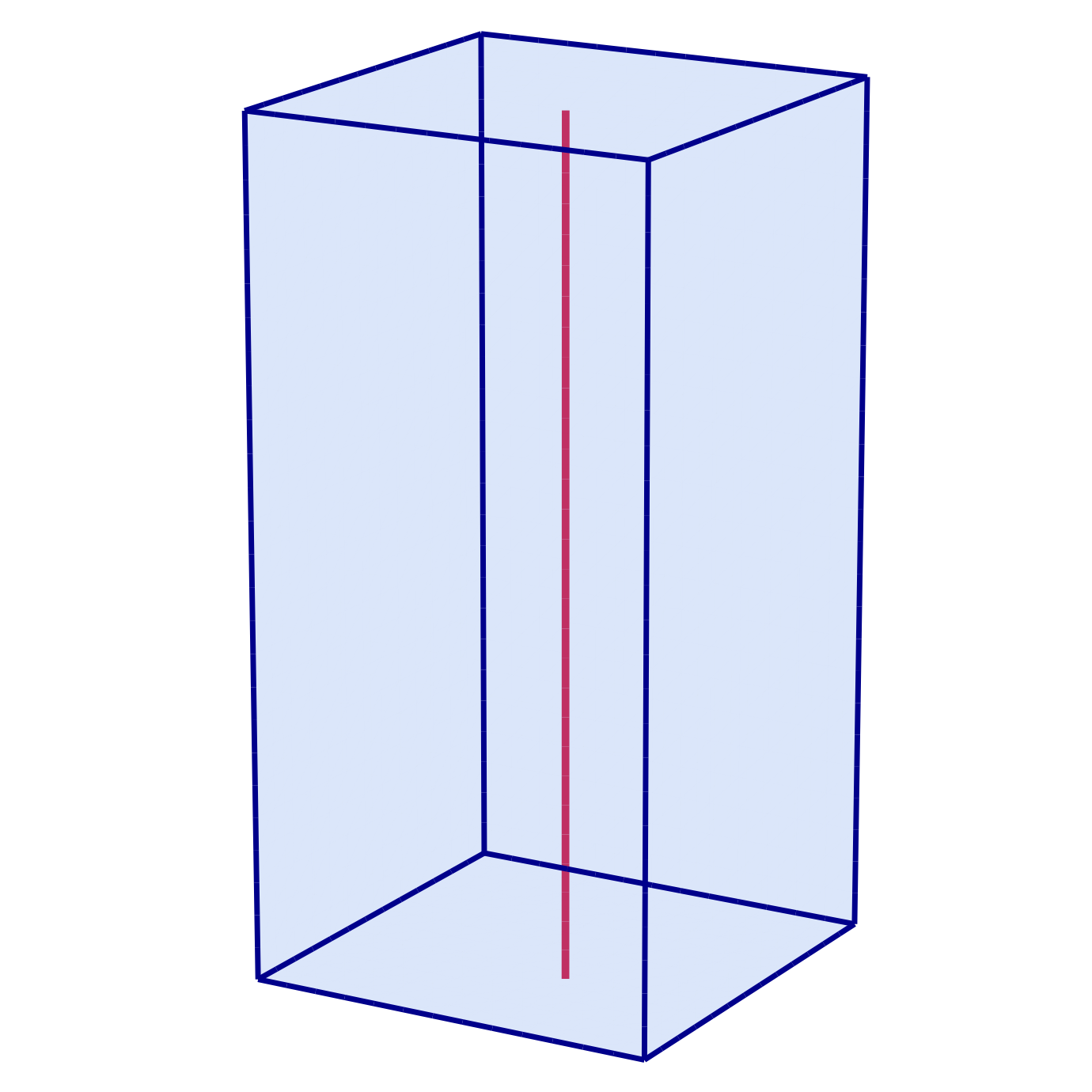}
  \includegraphics[width=0.31\textwidth]{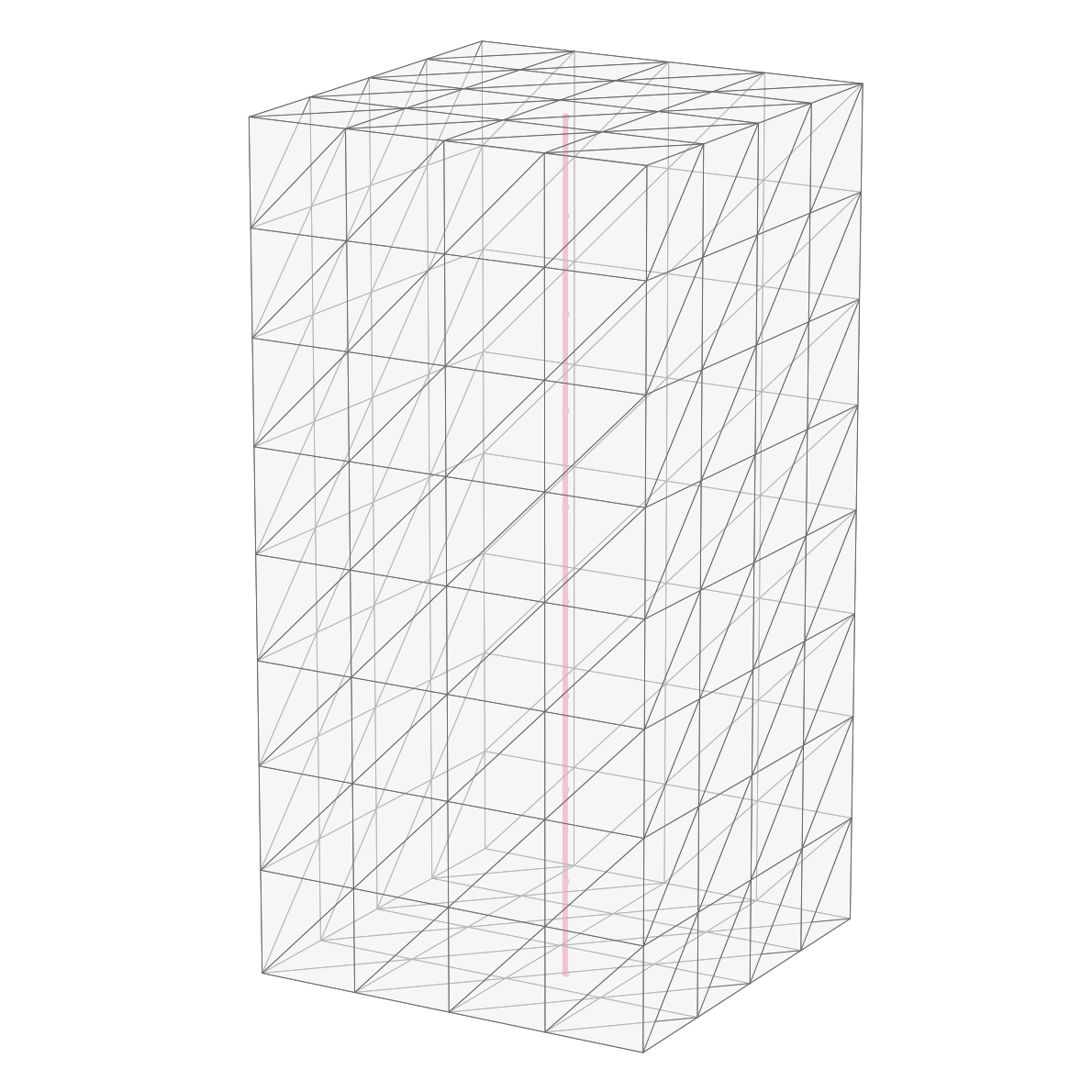}
  \includegraphics[width=0.31\textwidth]{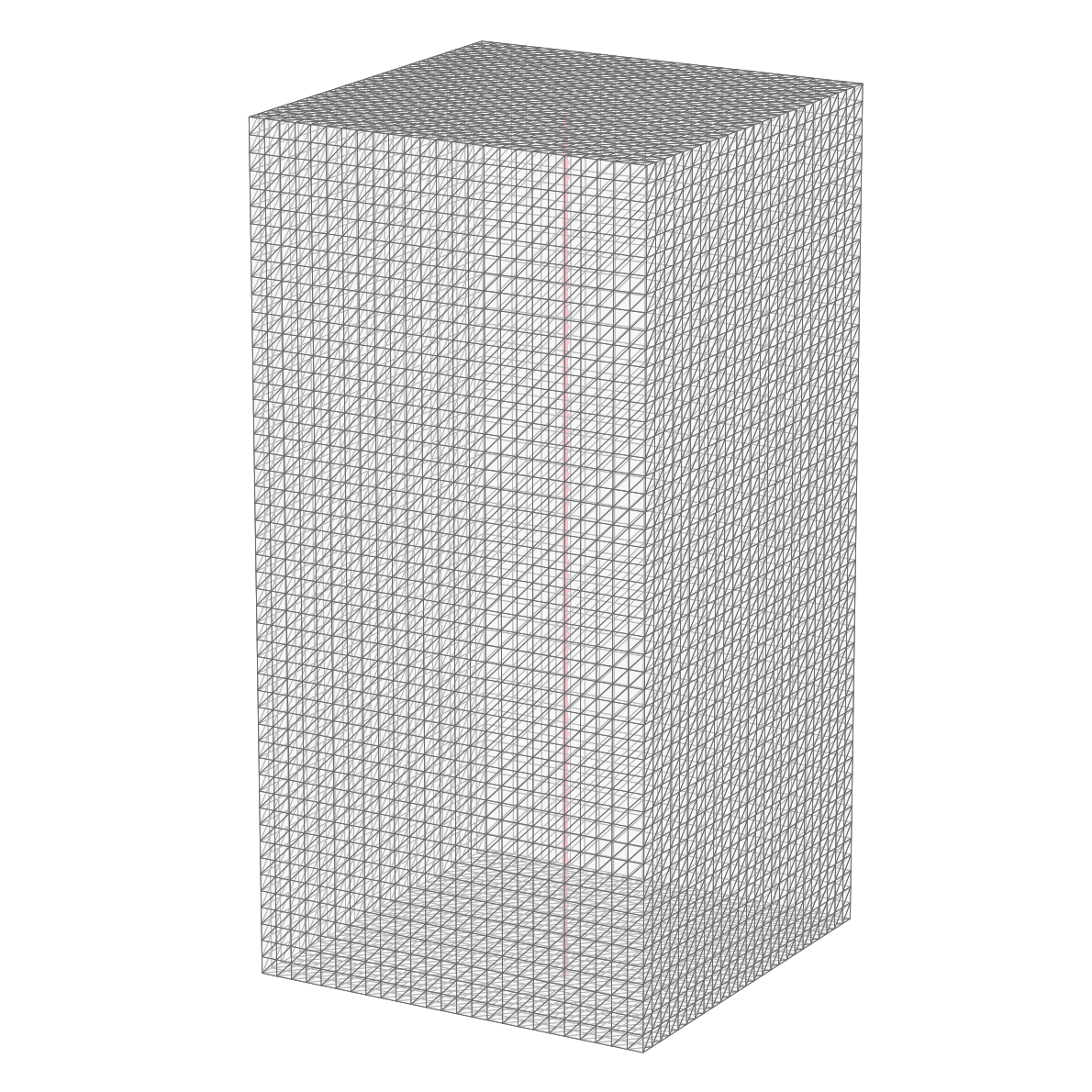}
  \caption{Left: brick with the off-axis embedded fiber, running
    its full height. Right: coarsest and finest solid meshes used in the refinement study, fiber
    mesh overlaid.}
  \label{fig-quadratic-geometry} 
\end{figure}

We embed a straight fiber, parallel to one edge but off-axis, in a conductive brick
(see Figure~\ref{fig-quadratic-geometry})
$[-\tfrac12,\tfrac12]^2\times[-1,1]$, at $(x_0,y_0)=(0.1,-0.15)$. The solid conductivity
is $\kappa=1$, and it is subject to a volumetric heat source $q=4$ and constant Neumann flux on all
six faces. No Dirichlet conditions are imposed on the solid, so its temperature is fixed
only through the tie to the fiber, whose two ends are prescribed.
The fiber has cross section $A=\pi\,0.05^2$ and conductivity per unit length
$K=10A\kappa$. With these loads, the coupled problem admits the exact solution
\begin{equation}
  u(x,y,z) = 1+z-1.5x^2-0.5y^2\ , \quad
  \theta(z) = z + 0.97291667\ , \quad \lambda= 0\ ,
  \label{eq-quadratic-exact}
\end{equation}
where $u$ is quadratic in the $(x,y)$ coordinates (and hence, unlike the linear field of
Section~\ref{subs-patch-example}, not exactly representable on linear tetrahedra) and $\theta$
is the corresponding fiber field, linear in $z$. The
multiplier vanishes because both the value gap and the axial-gradient gap between $u$ and
$\theta$ vanish identically. This tests the coupled discretization at every mesh level: the quadratic solid field carries a real discretization error that the discrete multiplier
$\lambda_h\neq0$ must balance, and its convergence to the exact $\lambda\equiv0$ is itself a verification measure.

The additive constant in $\theta$ deserves a word, since it is where the caveat of
Section~\ref{subs-implementation} becomes concrete. The constraint ties $\theta$ not to the value of
$u$ on the axis, but to the transverse average of $u$ that the coupling actually samples. For the
quadratic field at hand the axial value would give $z+0.97375$ and a rule integrating exactly over
the disk would give $z+0.9725$, whereas the two-diameter rule of Appendix~\ref{appendix}, which
assigns the second moment $R^2/6$ to each transverse direction, gives the constant shown
in~Eq.\eqref{eq-quadratic-exact}. The three differ by $O(R^2)$, i.e.\ by the square of the slenderness,
and any of them is an equally legitimate statement of the model; what matters here is that
\eqref{eq-quadratic-exact} is the exact solution of the coupled problem \emph{as discretized},
coupling operator and its quadrature included. That is precisely what this example sets out to
verify, and it is the only version of the statement against which a convergence rate is meaningful.

We solve the problem on four levels of uniform solid mesh refinement, refining the fiber mesh
alongside the solid mesh at every level, while ensuring $h_{\mathcal{C}}< h_{\Omega}$. Since the
exact solution is polynomial and the elements are affine, every error integral below is evaluated
exactly in closed form, so the reported errors carry no quadrature error of their own. As in
Section~\ref{subs-patch-example}, the multipliers at the two fiber ends, where the fiber temperature
is prescribed, are removed so that $W_h$ is exactly the discrete fiber temperature space.

\begin{table}[htbp]
  \centering
  \caption{Convergence under uniform refinement. The reported rate is the least-squares fit across
    all levels.
$\|e_u\|_{L^2}$: $L^2$ norm of error in $u$;
$|e_u|_{H^1}$: $H^1$ seminorm of error in $u$;
$\|e_\theta\|_{L^2}$: $L^2$ norm of error in $\theta$;
$\max|\lambda_h|$: largest nodal multiplier, whose exact value is zero.
  }
  \label{tab:convergence}
  \begin{tabular}{rrrrr}
    \toprule
    $h$ & $\|e_u\|_{L^2}$ & $|e_u|_{H^1}$ & $\|e_\theta\|_{L^2}$ & $\max|\lambda_h|$ \\
    \midrule
    0.315 & $1.484\cdot10^{-2}$ & $1.630\cdot10^{-1}$ & $3.840\cdot10^{-4}$ & $3.862\cdot10^{-1}$ \\
    0.18 & $4.141\cdot10^{-3}$ & $9.502\cdot10^{-2}$ & $1.091\cdot10^{-4}$ & $1.468\cdot10^{-1}$ \\
    0.09693 & $1.229\cdot10^{-3}$ & $5.162\cdot10^{-2}$ & $3.026\cdot10^{-5}$ & $3.289\cdot10^{-2}$ \\
    0.0504 & $4.578\cdot10^{-4}$ & $2.693\cdot10^{-2}$ & $7.996\cdot10^{-6}$ & $6.862\cdot10^{-3}$ \\
    \midrule
    Fitted slope & 1.90 & 0.98 & 2.11 & 2.23 \\
    Theoretical & 2 & 1 & 2 & --- \\
    \bottomrule
  \end{tabular}
\end{table}

\begin{figure}[htbp]
  \centering
  \includegraphics[width=0.75\textwidth]{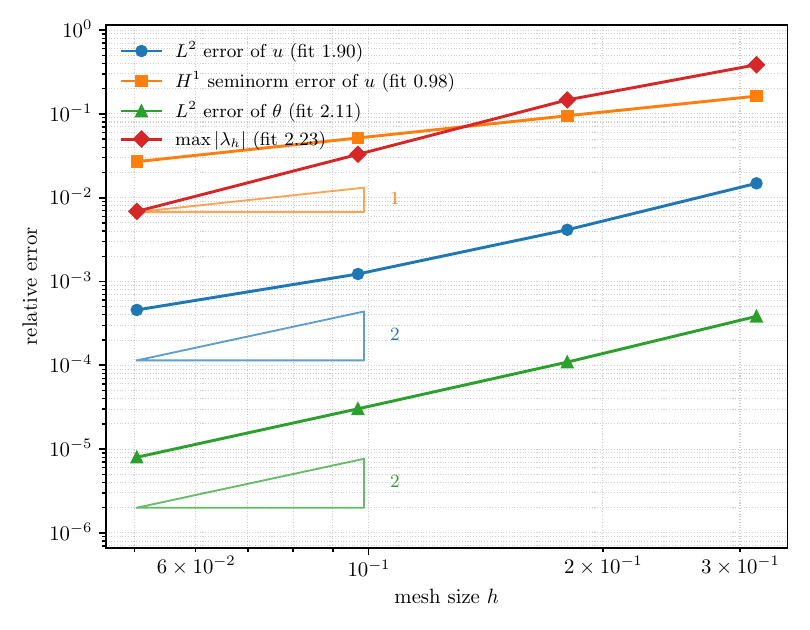}
  \caption{Convergence of the solid $L^2$ and $H^1$-seminorm
    errors, the fiber $L^2$ error, and $\max|\lambda_h|$, under uniform refinement. Triangles
    show the theoretical slopes $2$ ($L^2$) and $1$ ($H^1$).}
  \label{fig-quadratic-convergence} 
\end{figure}

\begin{figure}[htbp]
  \centering
  \includegraphics[height=0.4\textwidth]{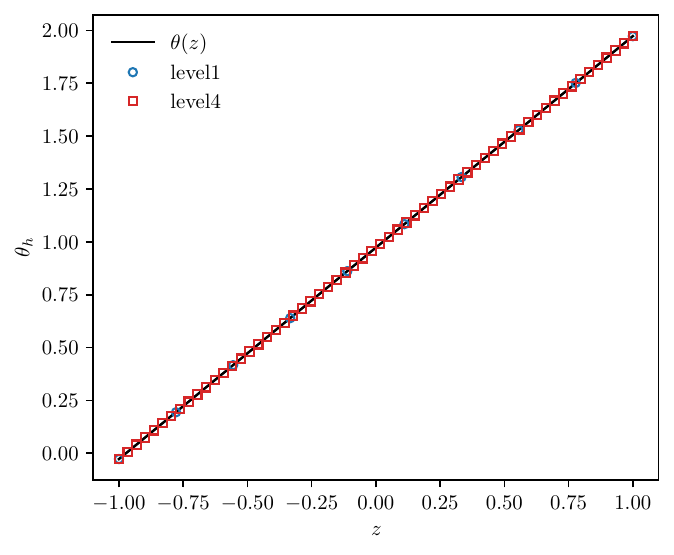}
  \hfill
  \includegraphics[height=0.4\textwidth]{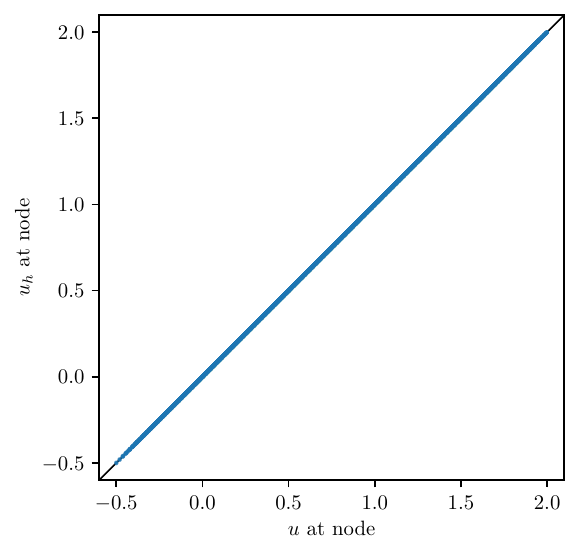}
  \caption{Left: computed fiber temperature at the coarsest and
    finest levels against the exact line $\theta(z)$. Right: parity plot of every solid nodal
    temperature at the finest level against the exact field.}
  \label{fig-quadratic-solution} 
\end{figure}

Table~\ref{tab:convergence} and Figure~\ref{fig-quadratic-convergence} show that the formulation
converges with the correct rates. The solid temperature converges at a fitted rate $1.90$ in the
$L^2$ norm and at a rate of $0.98$ in the $H^1$ seminorm, the fiber temperature converges at rate
$2.11$, and the discrete multiplier converges to its exact zero value at rate $2.23$.
Figure~\ref{fig-quadratic-solution} shows this agreement directly: the fiber temperature at the
coarsest and finest levels against the exact line $\theta(z)$, and a parity plot of every solid
nodal value against the exact field at the finest level.

\subsection{A helix with a highly conductive core}
\label{subs-helix-example}

The two previous examples verified the formulation on geometries where the
fiber is straight (Sections~\ref{subs-patch-example}
and~\ref{subs-quadratic-example}): the coupling samples along the fiber's
length always see a solid field that varies only with a single axial
coordinate. The next example checks that nothing about the discretization
depends on that simplification, by embedding the fiber --- and, in this
case, the solid itself --- along a curved, non-planar path: a
coil of three complete turns.

\begin{figure}[htbp]
  \centering
  \includegraphics[width=0.28\textwidth]{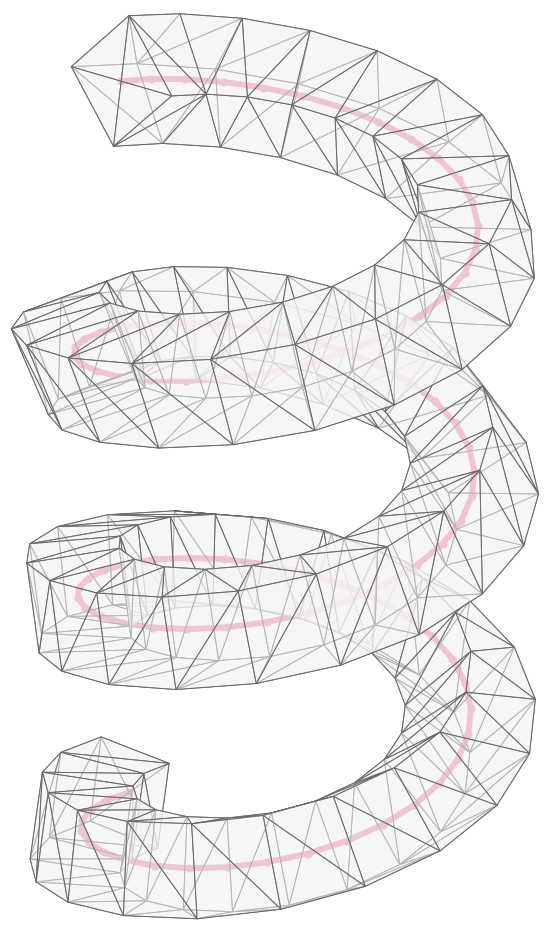}\hfill
  \includegraphics[width=0.28\textwidth]{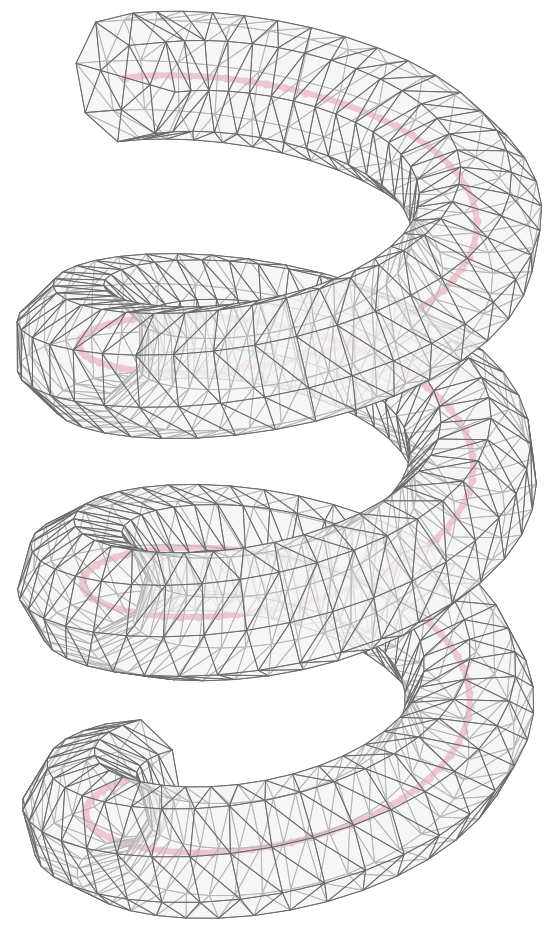}\hfill
  \includegraphics[width=0.28\textwidth]{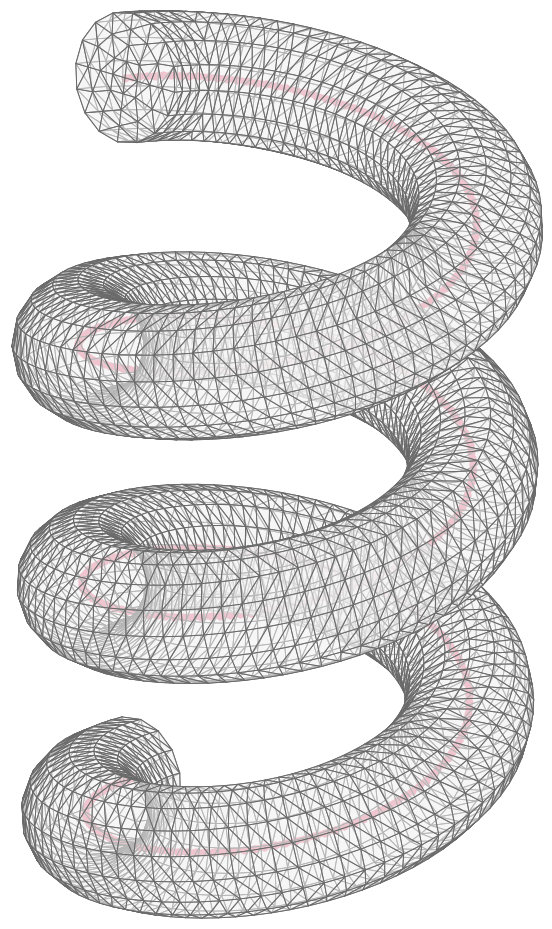}
  \caption{The three meshes used in the convergence study.
    Embedded fiber is shown in color. Mesh in reference solution is not shown.}
  \label{fig-helix-geometry}
\end{figure}

The solid is a helicoidal cylinder: its centerline has radius $3$ about the
$z$-axis, completes three full turns with pitch $4$, and its circular cross section has
radius~$1$. The fiber runs along the same centerline, with a 
circular cross section of radius $R=0.1$.
The solid mesh consists of tetrahedra and the fiber mesh employs
two-node linear elements (see Figure~\ref{fig-helix-geometry}).

The solid conductivity is $\kappa=1$ and the fiber conductivity per unit length
is $K=1$; the fiber's two ends are held at $\theta=1$ and
$\theta=0$, and the solid has no Dirichlet data, so its temperature
comes from the coupling with the wire. Unlike the example of Section~\ref{subs-quadratic-example},
there is no closed-form solution of the coupled problem on a helicoidal domain, so
verification here is by \emph{self-convergence}: the solid field $u_h$ at
each of three mesh levels is compared against a much finer reference mesh
(eight times as many solid tetrahedra as the finest of the three levels)
rather than against an exact solution. Because both the fiber's circular
cross section and its centerline are exact analytic images at every
level, a coarser level's solid nodes lie exactly on the same domain as the
reference mesh, so the reference field can be linearly interpolated at
each coarse node and compared directly. At every level the fiber mesh is
refined together with the solid mesh, keeping the fiber element length
below half the solid element size.

\begin{table}[htbp]
  \centering
  \caption{Convergence under uniform refinement. Rates are computed between consecutive levels.}
  \label{tab:convergence-helix}
  \begin{tabular}{rrcrc}
    \toprule
    \multirow{2}{*}{$h$} & \multicolumn{2}{c}{$L^2$ error of $u$ vs. finest-level reference} & \multicolumn{2}{c}{$\max|\lambda_h|$ (diagnostic)} \\
    \cmidrule(lr){2-3} \cmidrule(lr){4-5}
    & error & rate & error & rate \\
    \midrule
    1.262 & $1.913\cdot10^{-2}$ & --- & $4.634\cdot10^{0}$ & --- \\
    0.7356 & $7.874\cdot10^{-3}$ & 1.65 & $5.077\cdot10^{0}$ & -0.17 \\
    0.3983 & $2.479\cdot10^{-3}$ & 1.88 & $3.908\cdot10^{0}$ & 0.43 \\
    \midrule
    \multicolumn{1}{r}{fitted slope} & \multicolumn{2}{c}{1.77} & \multicolumn{2}{c}{0.15} \\
    \multicolumn{1}{r}{theoretical} & \multicolumn{2}{c}{2} & \multicolumn{2}{c}{---} \\
    \bottomrule
  \end{tabular}
\end{table}

\begin{figure}[htbp]
  \centering
  \includegraphics[width=0.65\textwidth]{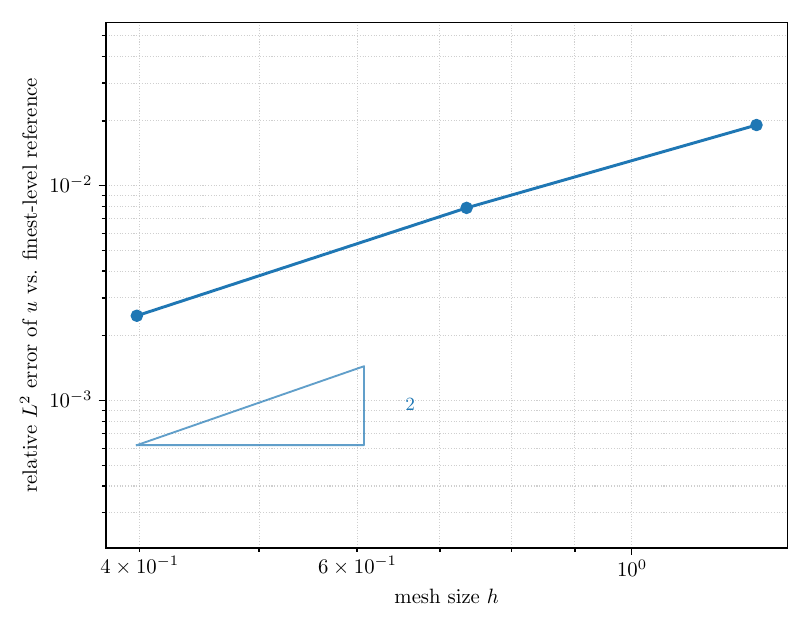}
  \caption{\label{fig-helix-convergence} Convergence of the solid
    field $u$ against the solution on the finest reference mesh.}
\end{figure}

\begin{figure}[htbp]
  \centering
  \includegraphics[width=0.27\textwidth]{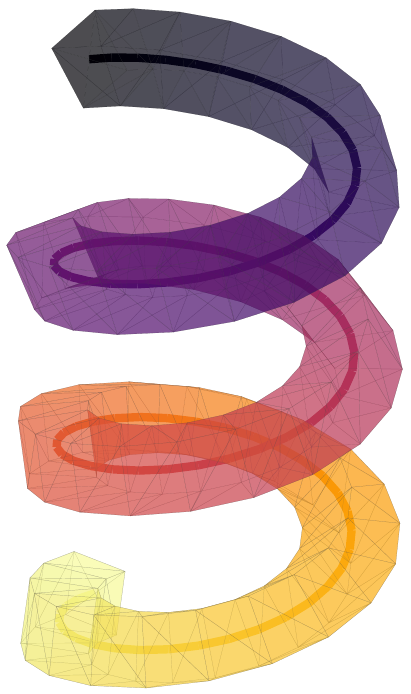}\hfill
  \includegraphics[width=0.27\textwidth]{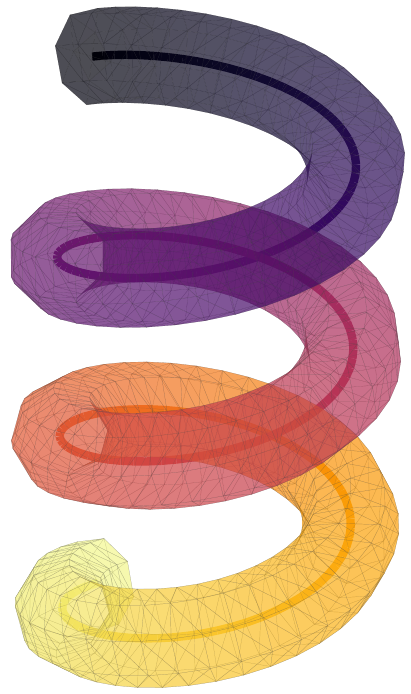}\hfill
  \includegraphics[width=0.27\textwidth]{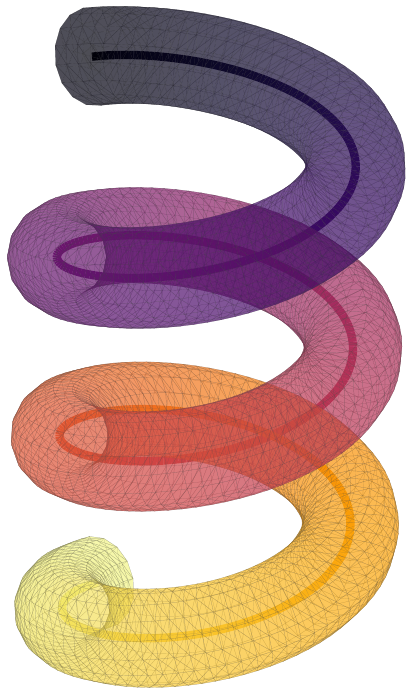}

\vspace{0.5cm} \includegraphics[width=0.6\textwidth]{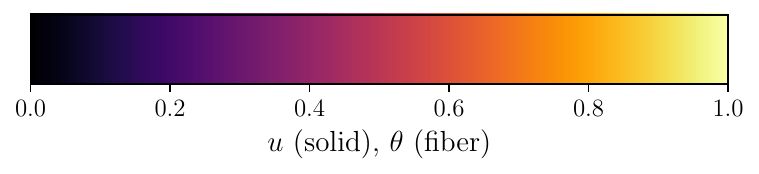}
  \caption{\label{fig-helix-contour} Computed temperature field at mesh
    levels 1, 2, and 3 (coarsest to finest): the solid field $u$ on the
    boundary surface and the fiber field $\theta$ at its core, on a shared
    color scale (below).}
\end{figure}

Table~\ref{tab:convergence-helix} and Figure~\ref{fig-helix-convergence}
show that the solid field converges at a fitted rate of $1.77$ in the
$L^2$ norm, approaching the theoretical rate $2$ for piecewise-linear
elements (the level-to-level rates are $1.65$ and $1.88$). The multiplier max value $\max|\lambda_h|$, also reported in
the table, does not converge to zero as it did in
Section~\ref{subs-quadratic-example}; this is expected rather than a
defect, since here the fiber field is not manufactured to make the
coupling's constraint gap vanish, so there is no reason for $\lambda_h$ to
do so either.
Figure~\ref{fig-helix-contour} shows the computed solid and fiber temperature fields at the
three mesh levels on a shared color scale.

\subsection{Conductive networks}
In this final example, we illustrate the generality of the proposed methodology by studying heat
transport in complex, branched, tree-like, networks embedded in conductive environments. Naturally,
these examples do not have a closed form solution so we will not attempt to study the convergence of
their solutions but rather focus on the possibilities of our embedding method. Moreover, we consider
now \emph{transient} thermal solutions to illustrate that the proposed ideas translate to these
situations. In fact, the linking terms of the coupled formulation are identical for stationary and
transient problem and the only difference between the previous examples and the current one is that
the formulation of the matrix thermal behavior now includes transient terms. We ignore the transient
contribution of the fiber.

We consider {four} different transport networks generated pseudo-randomly with a
space-colonization (attraction-point) growth algorithm
(\cite{runions2005ki,runions2007za}), each rooted near the bottom of a shared
$10\times10\times14$ box domain and differing only in their growth parameters: the strength of the
tropism toward the growth direction, the angular noise at each branching event, and the number,
radius, and reach of the attraction points that guide the growth (Table~\ref{tab-tree-params}).
All four share the same host solid (conductivity $1$, capacity $2$) and the same coupling
(cross section $A=0.05$, fiber conductivity $K=50$). Each root is held at $h=100$ from
$t=0$, all box walls are insulated, and the transient problem is integrated with the backward
Euler scheme to $t=100$ and fixed time step size $\Delta t=5$; because a pure-Neumann problem with a single Dirichlet point has no
steady solution other than a spatially uniform one, each network is instead snapshotted well
before equilibrium, at $t=10,30,60,100$, so that its shape remains visible as a thermal halo.

\begin{table}[htbp]
  \centering
  \caption{Growth parameters distinguishing the four networks.}
  \label{tab-tree-params} 
  \small
  \begin{tabular}{lrrrr}
    \toprule
    Network & tropism & branching noise & attr.\ points & attr.\ radius \\
    \midrule
    Column    & $0.85$ & $0.08$ & $150$ & $2.5$ \\
    Bush      & $0.15$ & $0.50$ & $500$ & $2.0$ \\
    Canopy    & $0.60$ & $0.15$ & $350$ & $15.0$ \\
    Windswept & $0.50$ & $0.25$ & $180$ & $2.5$ \\
    \bottomrule
  \end{tabular}
\end{table}

\begin{figure}[p]
  \centering
  \includegraphics[width=\textwidth]{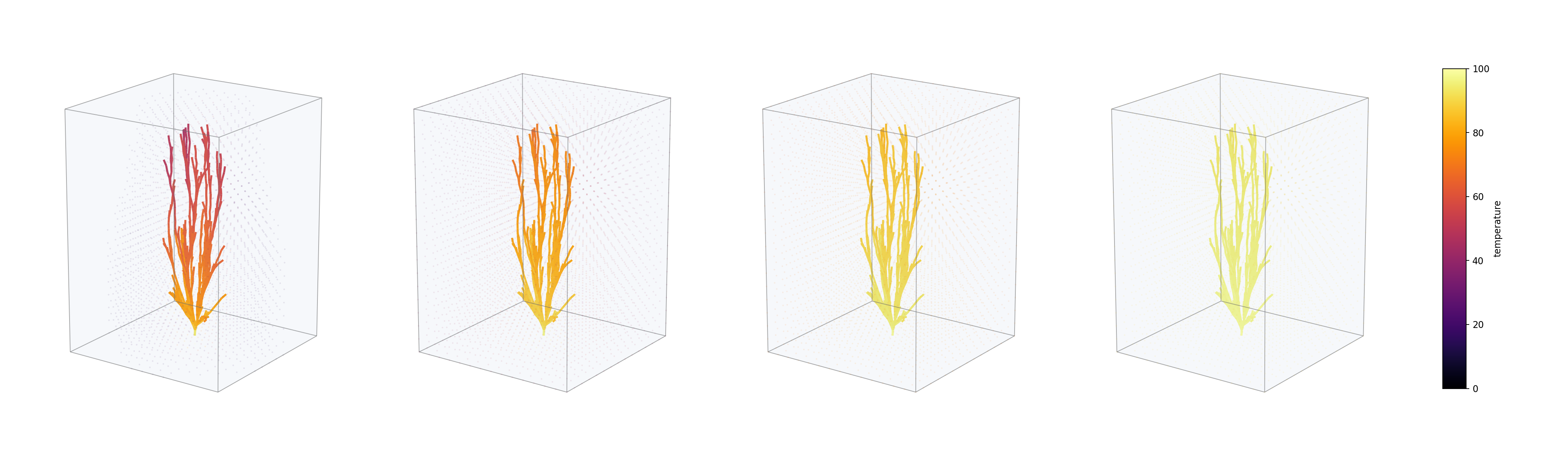}
  \caption{Column network: $t=10,30,60,100$ (left to right).}
  \label{fig-tree-column} 
\end{figure}

\begin{figure}[htbp]
  \centering \includegraphics[width=\textwidth]{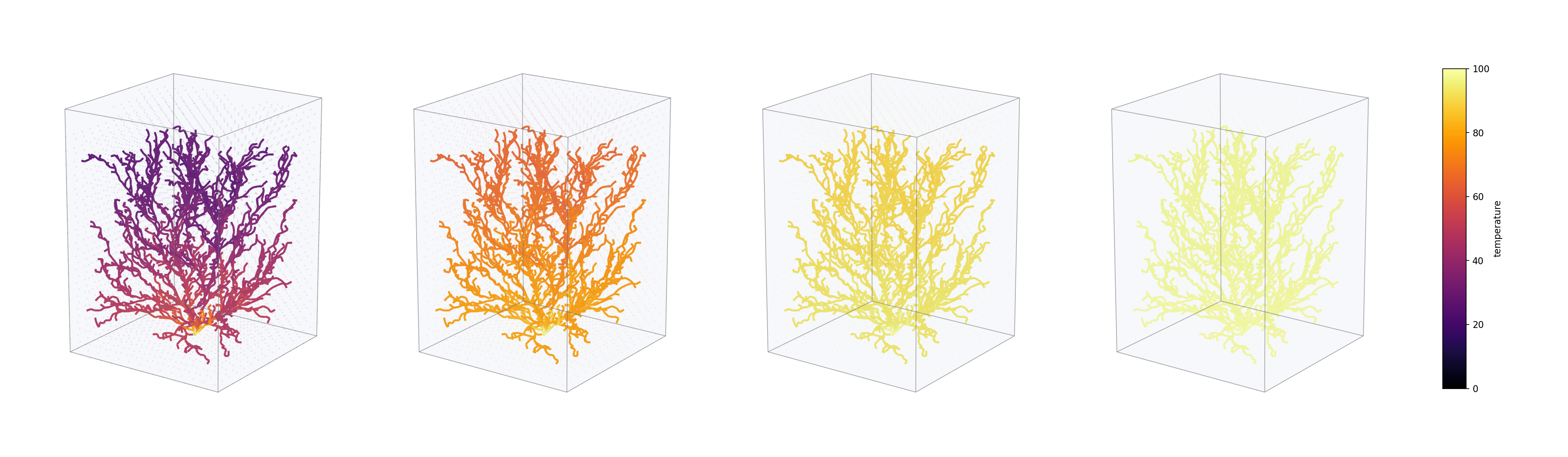}
  \caption{Bush network: $t=10,30,60,100$ (left to right).}
  \label{fig-tree-bush} 
\end{figure}

\begin{figure}[htbp]
  \centering
  \includegraphics[width=\textwidth]{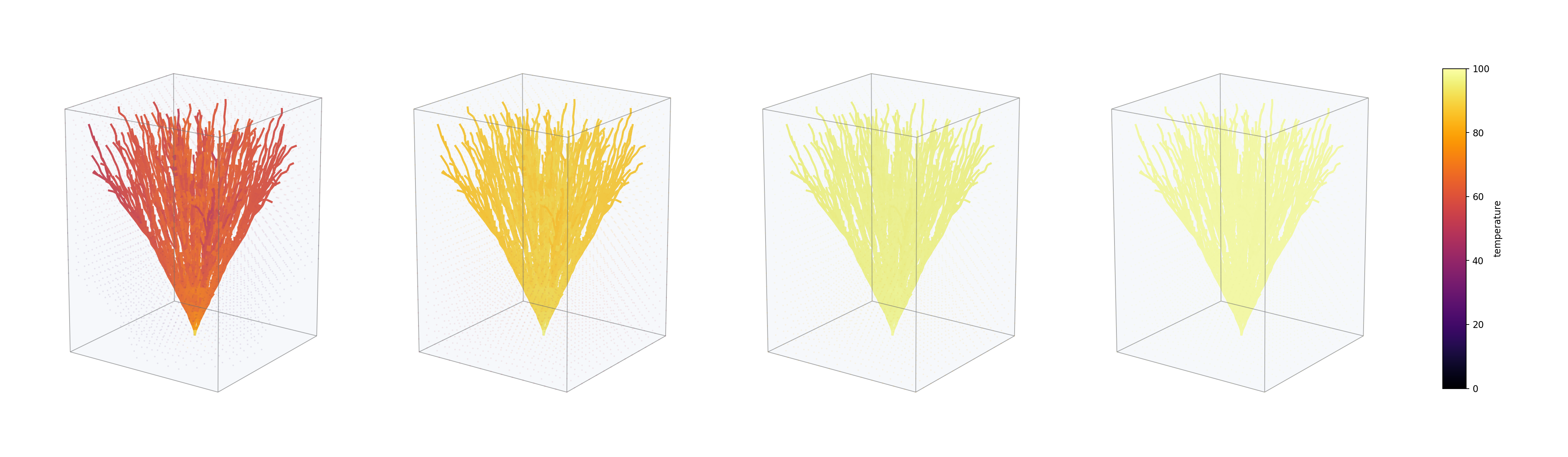}
  \caption{\label{fig-tree-canopy} Canopy network: $t=10,30,60,100$ (left to right).}
\end{figure}

\begin{figure}[htbp]
  \centering
  \includegraphics[width=\textwidth]{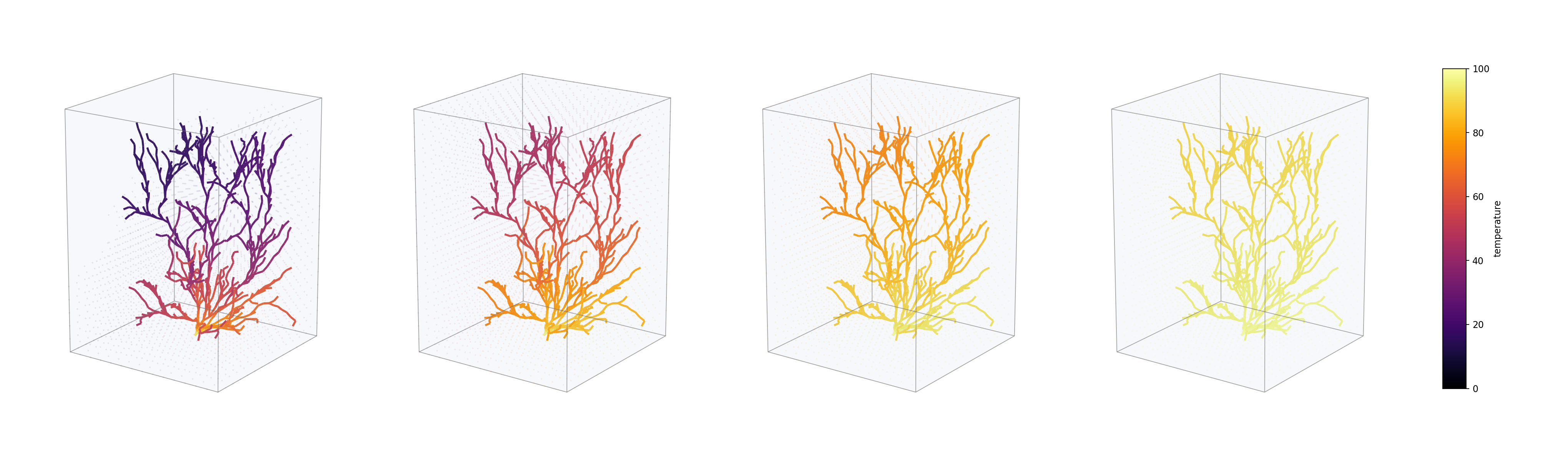}
  \caption{\label{fig-tree-windswept} Windswept network: $t=10,30,60,100$ (left to right).}
\end{figure}

High tropism with low branching noise (Column) gives a mostly straight, sparsely-branched trunk;
low tropism with high noise and a dense attraction-point cloud (Bush) gives a highly branched,
space-filling shape; a wide attraction radius (Canopy) lets branches reach far from the trunk
before terminating, producing a spreading crown; and a tilted growth direction (Windswept, not
shown in the table) biases the whole structure off-axis. Beyond their visual variety, the four
networks also stress the coupling's robustness on genuinely unstructured, non-conforming geometry:
between $10344$ and $88608$ transversely-sampled coupling points per network, of which all but
$15$ --- all in the windswept case, whose tilted growth direction lets branches approach the box's
corners more closely than the axis-aligned trees --- land inside the solid mesh, and every one of
the $20$ backward-Euler steps converges in a single Newton iteration in all four cases, with the
system's energy decaying monotonically toward the eventual uniform-temperature equilibrium
throughout. Figures~\ref{fig-tree-column}--\ref{fig-tree-windswept}
show, for each network, the four snapshots side by side on a common color scale. In every case
the thermal halo starts concentrated near the root and spreads outward as heat continues to be
supplied, tracing out the network's own branching structure before blurring into a more diffuse
cloud by $t=100$; the shape of the halo visibly reflects each network's own morphology, confirming
that the coupling correctly transports heat from an arbitrarily branched one-dimensional network
into the surrounding three-dimensional solid without requiring the two meshes to be compatible.

\section{Conclusions}
\label{sec-conclusions}
We have presented in this article a boundary value problem that describes the coupled
behavior of a slender diffusive body embedded inside a three-dimensional body that has a similar
diffusive behavior. This problem is based on our previous work on embedded structures which, in
turn, uses ideas of the Arlequin method. The proposed formulation is fairly general and can be shown
to lead to a well-posed saddle point optimization problem whose solution corresponds to the
coupled diffusive fields on the matrix and fiber.

One of the most attractive features of our proposed model is that it can be easily discretized
using finite elements, and the latter are unconditionally stable and convergent. By using a special
integration rule in the interface region, the meshes in the fiber and the matrix may be selected
independently, thus simplifying the meshing steps of an analysis. In addition, even though the
method has been presented in the context of a single embedded thermal conductor, its applicability
is much larger: it can be used, almost without modifications to other diffusive problems, with
multiple independent or connected fibers, and also with partially embedded ones.

The numerical examples that we have shown illustrate the properties that we have discussed. The
formulation is computationally efficient and, especially, robust.

\section{Acknowledgements}
Both authors acknowledge the funding received from the Spanish Ministry of Science, Innovation, and
Universities under grant PID2025-174275NB-I00

\section*{Declaration of generative AI and AI-assisted technologies in the manuscript preparation process}
During the preparation of this work, IR used Claude v2.1 to write scripts for automating the numerical examples. After using this tool, the author reviewed and edited the content as needed and takes full responsibility for the content of the published article.

\bibliographystyle{unsrt}
\bibliography{thermal}

\appendix
\section{Implementation details}
\label{appendix}
We provide in this appendix a compact pseudo code that describes all
the steps in the finite element implementation of the coupling terms of the proposed method.
The terms that come from the discretization of the Poisson problem in the matrix and
the fiber are standard and not described here.

\subsection{Basics}
A one-dimensional thermal bar (a fiber) is embedded inside a three-dimensional
solid mesh and coupled to it through a Lagrange-multiplier constraint. The bar has
circular cross-section area $A$ with radius $R=\sqrt{A/\pi}$. The finite element
implementation of the coupling terms has two stages,
corresponding to the two algorithms below. For simplicity, we show only the details
when the solid mesh consists of tetrahedra and the fiber uses two-node linear elements.

\begin{itemize}
  \item \textbf{Element creation} (Algorithm~\ref{alg:create}). For every
        point sampled on the cylindrical region $\mathcal{S}$, the solid
        element containing it is located, and a discrete coupling element
        is built connecting the four solid nodes, the two
        bar nodes, and two auxiliary Lagrange-multiplier nodes.
        
  \item \textbf{Element evaluation} (Algorithm~\ref{alg:evaluate}). Each
        link element evaluates, at its single sample point, an
        energy that enforces the zero $H^1$ gap between matrix and fiber
        temperatures. Its residual and tangent follow by
        differentiation.
\end{itemize}

\subsection{Notation}

\begin{center}
\begin{tabular}{@{}ll@{}}
\toprule
Symbol & Meaning \\
\midrule
$\mbs{x}$  & coordinates of a sample point in the region $\mathcal{S}$ \\
$\mbs{\eta}$ & barycentric coordinates of $\mbs{x}$ in its host tetrahedron \\
$\zeta\in[-1,1]$        & parametric coordinate of the sample point along the bar element \\
$\omega$               & integration weights of the sample point $\omega_{ij}=W_{i}w_j$ \\
$N_a$     & linear shape functions of the host tetrahedron, $a=1,\dots,4$ \\
$M_b$     & linear shape functions of the bar element, $b=1,2$ \\
$\mbs{t}$          & unit tangent of the bar element \\
$u_h, \theta_h$   & temperature fields interpolated on the solid and on the bar \\
$\lambda_h$       & Lagrange multiplier field \\
$g$                    & temperature gap, $g = u_h - \theta_h$ (see the note on signs below) \\
\bottomrule
\end{tabular}
\end{center}

\subsection{Algorithm 1: Creation of link elements}

The bar is discretized into 2-node elements; each is sampled at two Gauss
points along its length, and at each of those, on a small set of points
distributed over a disk. Every such sample point that falls inside the
solid mesh produces one coupling element.

\begin{algorithm}[H]
\caption{Creation of embedded thermal-bar coupling elements}
\label{alg:create}
\begin{algorithmic}[1]
\Procedure{CreateThermalBarLinks}{solid mesh, bar mesh, area $A$, subdivisions $n_h$}
  \State $R \equiv \ell \gets \sqrt{A/\pi}$
  \ForAll{bar elements $b_e$}
    \State compute the bar's unit tangent and two unit vectors spanning its cross section
    \ForAll{axial Gauss points $\xi_{\text{ip}}$ on $b_e$ (2-point rule)}
      \State $\mbs{x}_0 \gets$ physical position of $\xi_{\text{ip}}$ on the bar centerline
      \ForAll{transverse directions $\alpha$ (perpendicular diameters of the cross section)}
        \ForAll{transverse sample points $\xi_{\text{ih}}$, $\xi_{\text{ih}}\in[-1,1]$, $n_h+1$ points}
          \State $\mbs{x} \gets \mbs{x}_0 + \xi_{\text{ih}} \, R \, \text{(direction }\alpha\text{)}$
          \State $\omega \gets A \times (\text{axial weight} \times \text{axial jacobian}) \times \tfrac14 \times (\text{transverse weight})$
          \State \Comment{find the solid element $e$ in which point $\mbs{x}$ falls}
          \State $e, \mbs{\eta} \gets$ \Call{LocateHostElement}{$\mbs{x}$, solid mesh}
          \If{no host element found}
            \State issue warning and skip this sample point
          \Else
            \State record: host solid element $e$, its barycentric coords $\mbs{\eta}$,
                   the bar element $b_e$, the parametric point $\xi_{\text{ip}}$, and the weight $\omega$
            \State mark $b_e$ as ``embedded'' (its two nodes will need Lagrange multipliers)
          \EndIf
        \EndFor
      \EndFor
    \EndFor
  \EndFor
  \State \Comment{one Lagrange-multiplier node per embedded bar node, collocated with it}
  \ForAll{nodes $n$ of embedded bar elements}
    \If{no multiplier node created yet for $n$}
      \State create a new scalar-dof node at $n$'s position; associate it with $n$
    \EndIf
  \EndFor
  \ForAll{recorded samples $(e, \mbs{\eta}, b_e, \xi_{\text{ip}}, w)$}
    \State node list $\gets$ [4 nodes of $e$] $+$ [2 nodes of $b_e$] $+$ [2 multiplier nodes of $b_e$'s nodes]
    \State create a link element from the node list, $\mbs{\eta}$, $\omega$, $\xi_{\text{ip}}$
  \EndFor
  \State \Return set of created link elements
\EndProcedure
\end{algorithmic}
\end{algorithm}

\noindent Each resulting link element therefore has 8 nodes:
4 solid (temperature dof $u_h$), 2 bar (temperature dof $\theta_{h}$), and
2 Lagrange multiplier (dof $\lambda_h$), and carries as data the host
barycentric coordinates $\mbs{\eta}$, the bar parametric coordinate $\zeta$, and the
sample weight $\omega$.

\subsection{Algorithm 2: Lagrangian, residual, and tangent}

At its single sample point, the link element interpolates the solid temperature
$u_h$, the bar temperature $\theta_h$, and the multiplier $\lambda_h$, forms the temperature
gap $g$, and assembles the associated weak-form contributions.

One convention differs from the body of the article and is worth stating, since it changes the sign
of the computed multiplier. Eq.~\eqref{eq-lagrangian} writes the constraint with the gap
$\mathcal{L}\beta-v$, whereas the implementation below uses $g=u_h-\theta_h$, the opposite sign.
The two Lagrangians therefore differ by $\lambda\mapsto-\lambda$: the saddle point, the temperature
fields, and $|\lambda|$ are identical, but the sign of the multiplier reported by the code is
reversed with respect to the one in Section~\ref{subs-coupling}.

\begin{algorithm}[H]
\caption{Evaluation of Lagrangian, residual, and tangent of a link element}
\label{alg:evaluate}
\begin{algorithmic}[1]
\Procedure{Evaluate}{element with data $\mbs{\eta}, \zeta, \omega, \ell$}
  \State \Comment{shape functions at the fixed sample point}
  \State $N_a \gets$ barycentric coordinates $\mbs{\eta}$, $a=1,\dots,4$
  \State $M_b(\zeta) \gets \tfrac12(1\mp \zeta)$, $b=1,2$
  \State $M_b' \gets$ derivative of $M_b$ w.r.t.\ arclength
  \State $\nabla N_a \gets$ gradients of the solid shape functions
  \State $\mbs{t} \gets$ unit tangent of the bar element
  \State
  \State \Comment{interpolate fields}
  \State $u_{h} \gets \sum_a N_a \, u_a$, \quad $\theta_{h} \gets \sum_b M_b \, \theta_b$
  \State $\lambda_h \gets \sum_c M_c \, \lambda_c$
  \State $g \gets u_{h} - \theta_h$
  \State $\nabla u_h \gets \sum_a \nabla N_a \, u_a$, \quad
         $\theta_h' \gets \sum_b M_b' \, \theta_b$, \quad
         $\lambda_h' \gets \sum_c M_c' \, \lambda_c$
  \State $\nabla g \cdot \mbs{t} \gets \mbs{t}\cdot\nabla u_h - \theta_h'$
  \State
  \State \Comment{Contribution to Lagrangian}
  \State $L \gets \omega \Big[\, \lambda_h\, g \;+\; \ell^2 \, \lambda_h'\, (\nabla g \cdot \mbs{t}) \,\Big]$
  \State
  \State \Comment{residual: $\partial L/\partial(\cdot)$ for every dof of the element}
  \ForAll{solid node $a$}
    \State $R_a \mathrel{+}= \omega\Big[\lambda_h\,N_a + \ell^2\,\lambda_h'\,\mbs{t}\cdot\nabla N_a\Big]$
  \EndFor
  \ForAll{bar node $b$}
    \State $R_b \mathrel{+}= -\omega\Big[\lambda_h\,M_b + \ell^2\,\lambda_h'\,M_b'\Big]$
  \EndFor
  \ForAll{multiplier node $c$}
    \State $R_c \mathrel{+}= \omega\Big[M_c\, g + \ell^2\,M_c'\,(\nabla g\cdot \mbs{t})\Big]$
  \EndFor
  \State
  \State \Comment{tangent:}
  \ForAll{multiplier node $c$, solid node $a$}
    \State $K_{ca} \gets \omega\Big[M_c N_a + \ell^2\,M_c'\, \mbs{t}\cdot\nabla N_a\Big]$
    \State assemble $K_{ca}$ at $(\lambda_c, u_a)$ and, by symmetry, at $(u_a,\lambda_c)$
  \EndFor
  \ForAll{multiplier node $c$, bar node $b$}
    \State $K_{cb} \gets -\omega\Big[M_c M_b + \ell^2\,M_c'\; M_b'\Big]$
    \State assemble $K_{cb}$ at $(\lambda_c, \theta_b)$ and, by symmetry, at $(\theta_b,\lambda_c)$
  \EndFor
  \State \Return $L$, residual $R$, tangent $K$
\EndProcedure
\end{algorithmic}
\end{algorithm}

\end{document}